\documentstyle[11pt, amssymb]{article}

\begin{document}

\begin{flushleft}
{\Large {\bf Tau functions and KP solutions on families of algebraic curves}} 
\end{flushleft}

\begin{flushleft}
{\large {\bf Takashi Ichikawa}}  
\end{flushleft}

\noindent 
{\small Department of Mathematics, Faculty of Science and Engineering, 
Saga University, Saga 840-8502, Japan. 
E-mail: ichikawn@cc.saga-u.ac.jp} 
\vspace{2ex}

\noindent
{\bf Abstract:} 
Using abelian differentials and periods of the universal Mumford curve, 
we study the universal expression and asymptotic behavior of tau functions defined for 
stably degenerating families of algebraic curves with additional data. 
Furthermore, we apply this study to constructing solutions to the KP hierarchy 
which are expressed by nonarchimedean theta functions 
and by mixtures of quasi-periodic solutions with solitons 
containing real solutions defined on families of {\texttt M}-curves. 
\vspace{2ex}

\vspace{2ex}

\noindent
{\bf 1. Introduction} 
\vspace{2ex}

\noindent 
Tau functions were introduced by Sato, Jimbo, Miwa and others in the study of integrable systems. 
Especially, Riemann surfaces with additional data such as marked points and line bundles 
give rise to tau functions which play important roles in the theory of soliton equations 
(cf. \cite{Kr, SeW}) and of conformal field in theory in operator formalism 
(cf. \cite{Al-GGR, IMO, KNTY, V}). 
The aim of this paper is to study the tau functions defined for degenerating families of 
algebraic curves using their explicit formula in terms of theta functions 
from the viewpoint of formal arithmetic geometry. 
We have universal expressions of the tau functions and the associated solutions 
of the KP hierarchy (KP solutions for short) related with soliton solutions 
which are shown to be meaningful in the fields of complex geometry, 
real algebraic geometry, nonarchimedean geometry and tropical geometry. 

The organization of this paper is as follows. 

In Sections 2--5, 
as a main tool to showing our results, 
we review results given in \cite{I4, I5} on the universal Mumford curve 
with abelian differentials and periods, 
and its application to their variational formulas partially obtained in \cite{F, BM, HN}. 

In Section 6, 
we give universal and computable formulas of tau functions and 
the associated KP solutions for families of Schottky uniformized Riemann surfaces. 
By the arithmeticity of the formulas, 
extending results in \cite{I1}  
we can express tau functions for all Mumford curves 
over nonarchimedean complete valuation fields of characteristic $0$ 
in terms of the associated nonarchimedean theta functions (cf. \cite{BoL, FoRSS}). 

In Section 7, 
we study the asymptotic behavior of tau functions under degeneration of 
considering Riemann surfaces in order to extend the above conformal field theory 
for degenerated objects. 
From the variational formulas stated in Section 5, 
we can obtain modified tau functions for degenerate objects which give also KP solutions 
given by mixtures of quasi-periodic solutions and soliton ones. 
Similar results were obtained in \cite{N1}--\cite{N3} and \cite{BeEN} 
for hyperelliptic curves using the theory of Sato Grassmannian \cite{SS} and 
of sigma functions. 
Furthermore, generalizing results in \cite{Go, Mu3} 
we have modified tau functions expressed by rational functions 
which give rise to soliton solutions as a special case. 

In Section 8, 
we show that the universal Mumford curve give rise to degenerating families of 
{\texttt M}-curves, 
i.e., real algebraic curves with maximal number of components. 
Furthermore, we give a similar statement to results in \cite{AG} 
by showing that the tau functions associated with these {\texttt M}-curves 
are real, i.e., real-valued for real variables, 
and tending to soliton KP solutions. 
This result seems also related to theta functions of tropical curves and 
their limits which were discussed in \cite{AFMS}. 
\vspace{4ex}

\noindent
{\bf 2. Generalized Tate curve}  
\vspace{2ex}

\noindent
{\it 2.1. Schottky uniformization.} 
A Schottky group $\Gamma$ of rank $g$ is defined as a free group with generators 
$\gamma_{i} \in PGL_{2}({\mathbb C})$ $(1 \leq i \leq g)$ which map Jordan curves  
$C_{-i} \subset {\mathbb P}^{1}_{\mathbb C} = {\mathbb C} \cup \{ \infty \}$ 
to other Jordan curves $C_{i} \subset {\mathbb P}^{1}_{\mathbb C}$ 
with orientation reversed, 
where $C_{\pm 1},..., C_{\pm g}$ with their interiors are mutually disjoint. 
Each element $\gamma \in \Gamma - \{ 1 \}$ is conjugate 
to an element of $PGL_{2}({\mathbb C})$ sending $z$ to $\beta_{\gamma} z$ 
for some $\beta_{\gamma} \in {\mathbb C}^{\times}$ with $|\beta_{\gamma}| < 1$ 
which is called the {\it multiplier} of $\gamma$. 
Therefore, one has 
\begin{eqnarray*}
\frac{\gamma(z) - \alpha_{\gamma}}{z - \alpha_{\gamma}} = 
\beta_{\gamma} \frac{\gamma(z) - \alpha'_{\gamma}}{z - \alpha'_{\gamma}} 
\end{eqnarray*}
for some element $\alpha_{\gamma}, \alpha'_{\gamma}$ of 
${\mathbb P}^{1}_{\mathbb C}$ 
called the {\it attractive}, {\it repulsive} fixed points of $\gamma$ 
respectively. 
Then the discontinuity set 
$\Omega_{\Gamma} \subset {\mathbb P}^{1}_{\mathbb C}$ 
under the action of $\Gamma$ has a fundamental domain $D_{\Gamma}$ 
which is given by the complement of the union of $C_{\pm 1},..., C_{\pm g}$ with their interiors. 
The quotient space $R_{\Gamma} = \Omega_{\Gamma}/\Gamma$ 
is a (compact) Riemann surface of genus $g$ 
which is called {\it Schottky uniformized} by $\Gamma$ (cf. \cite{Sc}). 
Furthermore, by a result of Koebe, 
every Riemann surface of genus $g$ can be represented in this manner. 
\vspace{2ex}

\noindent
{\it 2.2. Generalized Tate curve.} 
A (marked) curve is called {\it degenerate} if it is a stable (marked) curve and 
the normalization of its irreducible components are all projective (marked) lines. 
Then the dual graph $\Delta = (V, E, T)$ of a stable marked curve is a collection  
of 3 finite sets $V$ of vertices, $E$ of edges, $T$ of tails 
and 2 boundary maps 
$$
b : T \rightarrow V, 
\ \ b : E \longrightarrow \left( V \cup \{ \mbox{unordered pairs of elements of $V$} \} \right) 
$$
such that the geometric realization of $\Delta$ is connected 
and that $\Delta$ is {\it stable}, namely its each vertex has at least $3$ branches. 
The number of elements of a finite set $X$ is denoted by $\sharp X$, 
and a (connected) stable graph $\Delta = (V, E, T)$ is called {\it of $(g, n)$-type} 
if ${\rm rank}_{\mathbb Z} H_{1}(\Delta, {\mathbb Z}) = g$, $\sharp T = n$. 
Then under fixing a bijection $\nu : T \stackrel{\sim}{\rightarrow} \{ 1, ... , n \}$, 
which we call a numbering of $T$, 
$\Delta = (V, E, T)$ becomes the dual graph of a degenerate curve of genus $g$ 
with $n$ points such that each tail $h \in T$ corresponds to the $\nu(h)$th marked point. 
In particular, a stable graph without tail is the dual graph of 
a degenerate (unmarked) curve by this correspondence. 
If $\Delta$ is trivalent, i.e. any vertex of $\Delta$ has just $3$ branches, 
then a degenerate $\sharp T$-marked curve with with dual graph $\Delta$ 
is maximally degenerate. 
An {\it orientation} of a stable graph $\Delta = (V, E, T)$ means 
giving an orientation of each $e \in E$. 
Under an orientation of $\Delta$, 
denote by $\pm E = \{ e, -e \ | \ e \in E \}$ the set of oriented edges, 
and by $v_{h}$ the terminal vertex of $h \in \pm E$ (resp. the boundary vertex of $h \in T)$. 
For each $h \in \pm E$, 
denote by let $| h | \in E$ be the edge $h$ without orientation. 

Let $\Delta = (V, E, T)$ be a stable graph. 
Fix an orientation of $\Delta$, 
and take a subset ${\mathcal E}$ of $\pm E \cup T$ 
whose complement ${\mathcal E}_{\infty}$ satisfies the condition that 
$$
\pm E \cap {\mathcal E}_{\infty} \cap \{ -h \ | \ h \in {\mathcal E}_{\infty} \} 
\ = \ 
\emptyset, 
$$ 
and that $v_{h} \neq v_{h'}$ for any distinct $h, h' \in {\mathcal E}_{\infty}$. 
We attach variables $x_{h}$ for $h \in {\mathcal E}$ and $y_{e} = y_{-e}$ for $e \in E$. 
Let $R_{\Delta}$ be the ${\mathbb Z}$-algebra generated by $x_{h}$ $(h \in {\mathcal E})$, 
$1/(x_{e} - x_{-e})$ $(e, -e \in {\mathcal E})$ and $1/(x_{h} - x_{h'})$ 
$(h, h' \in {\mathcal E}$ with $h \neq h'$ and $v_{h} = v_{h'})$, 
and let 
$$ 
A_{\Delta} \ = \ R_{\Delta} [[y_{e} \ (e \in E)]], \ \ 
B_{\Delta} \ = \ A_{\Delta} \left[ \prod_{e \in E} y_{e}^{-1} \right]. 
$$ 
According to \cite[Section 2]{I2}, 
we construct the universal Schottky group $\Gamma_{\Delta}$ 
associated with oriented $\Delta$ and ${\mathcal E}$ as follows. 
For $h \in \pm E$, put 
\begin{eqnarray*}
\phi_{h} 
& = & 
\left( \begin{array}{cc} x_{h} & x_{-h} \\ 1 & 1 \end{array} \right) 
\left( \begin{array}{cc} 1 & 0 \\ 0 & y_{h} \end{array} \right) 
\left( \begin{array}{cc} x_{h} & x_{-h} \\ 1 & 1 \end{array} \right)^{-1} 
\\ 
& = & 
\frac{1}{x_{h} - x_{-h}} 
\left\{ \left( \begin{array}{cc} x_{h} & - x_{h} x_{-h} \\ 1 & - x_{-h} \end{array} \right) - 
\left( \begin{array}{cc} x_{-h} & - x_{h} x_{-h} \\ 1 & - x_{h} \end{array} \right) y_{h} \right\}, 
\end{eqnarray*} 
where $x_{h}$ (resp. $x_{-h})$ means $\infty$ 
if $h$ (resp. $-h)$ belongs to ${\mathcal E}_{\infty}$. 
This gives an element of $PGL_{2}(B_{\Delta}) = GL_{2}(B_{\Delta})/B_{\Delta}^{\times}$ 
which we denote by the same symbol, 
and satisfies
$$
\frac{\phi_{h}(z) - x_{h}}{z - x_{h}} 
\ = \ 
y_{h} \frac{\phi_{h}(z) - x_{-h}}{z - x_{-h}} 
\ \ (z \in {\mathbb P}^{1}), \eqno(2.2.1)
$$
where $PGL_{2}$ acts on ${\mathbb P}^{1}$ by linear fractional transformation. 

For any reduced path $\rho = h(1) \cdot h(2) \cdots h(l)$ 
which is the product of oriented edges $h(1), ... ,h(l)$ such that $v_{h(i)} = v_{-h(i+1)}$, 
one can associate an element $\rho^{*}$ of $PGL_{2}(B_{\Delta})$ 
having reduced expression 
$\phi_{h(l)} \phi_{h(l-1)} \cdots \phi_{h(1)}$. 
Fix a base vertex $v_{b}$ of $V$, 
and consider the fundamental group 
$\pi_{1} (\Delta, v_{b})$ which is a free group 
of rank $g = {\rm rank}_{\mathbb Z} H_{1}(\Delta, {\mathbb Z})$. 
Then the correspondence $\rho \mapsto \rho^{*}$ 
gives an injective anti-homomorphism 
$\pi_{1} (\Delta, v_{b}) \rightarrow PGL_{2}(B_{\Delta})$ 
whose image is denoted by $\Gamma_{\Delta}$. 

It is shown in \cite[Section 3]{I2} and \cite[1.4]{I3} 
(see also \cite[Section 2]{IhN} when $\Delta$ is trivalent and has no loop) 
that for any stable graph $\Delta$, 
there exists a stable marked curve ${\mathcal C}_{\Delta}$ of genus $g$ over $A_{\Delta}$ 
which satisfies the following properties: 

\begin{itemize}

\item[(P1)] 
The closed fiber ${\mathcal C}_{\Delta} \otimes_{A_{\Delta}} R_{\Delta}$ of ${\mathcal C}_{\Delta}$ 
obtained by substituting $y_{e} = 0$ $(e \in E)$ 
becomes the degenerate marked curve over $R_{\Delta}$ with dual graph $\Delta$ which is 
obtained from the collection of $P_{v} := {\mathbb P}^{1}_{R_{\Delta}}$ $(v \in V)$ 
by identifying the points $x_{e} \in P_{v_{e}}$ and $x_{-e} \in P_{v_{-e}}$ ($e \in E$), 
where $x_{h}$ denotes $\infty$ if $h \in {\mathcal E}_{\infty}$. 

\item[(P2)] 
${\mathcal C}_{\Delta}$ gives rise to a universal deformation 
of degenerate marked curves with dual graph $\Delta$. 
More precisely, ${\mathcal C}_{\Delta}$ satisfies the following: 
For a noetherian and normal complete local ring $R$ with residue field $k$, 
let $C$ be a marked Mumford curve over $R$, 
namely a stable marked curve over $R$ with nonsingular generic fiber 
such that the closed fiber $C \otimes_{R} k$ is a degenerate marked curve 
with dual graph $\Delta$, in which all double points and marked points are $k$-rational. 
Then there exists a ring homomorphism $A_{\Delta} \rightarrow R$ 
giving ${\mathcal C}_{\Delta} \otimes_{A_{\Delta}} R \cong C$.  

\item[(P3)] 
${\mathcal C}_{\Delta} \otimes_{A_{\Delta}} B_{\Delta}$ is smooth over $B_{\Delta}$ 
and is Mumford uniformized (cf. \cite{Mu1}) by $\Gamma_{\Delta}$. 

\item[(P4)] 
Take $x_{h}$ $(h \in {\mathcal E})$ as complex numbers such that $x_{e} \neq x_{-e}$ 
and that $x_{h} \neq x_{h'}$ if $h \neq h'$ and $v_{h} = v_{h'}$, 
and take $y_{e}$ $(e \in E)$ as sufficiently small nonzero complex numbers. 
Then ${\mathcal C}_{\Delta}$ gives rise to a family of Riemann surfaces 
denoted by ${\mathcal R}_{\Delta}$ with marked points which are uniformized 
by Schottky groups over ${\mathbb C}$ obtained from $\Gamma_{\Delta}$. 

\end{itemize}
We review the construction of ${\mathcal C}_{\Delta}$ given in \cite[Theorem 3.5]{I2}. 
Let $T_{\Delta}$ be the tree obtained as the universal cover of $\Delta$, 
and denote by ${\mathcal P}_{T_{\Delta}}$ be the formal scheme 
as the union of ${\mathbb P}^{1}_{A_{\Delta}}$'s indexed by vertices of $T_{\Delta}$ 
under the $B_{\Delta}$-isomorphism by $\phi_{e}$ $(e \in E)$. 
Then it is shown in \cite[Theorem 3.5]{I2} that 
${\mathcal C}_{\Delta}$ is the formal scheme theoretic quotient of ${\mathcal P}_{T_{\Delta}}$ by $\Gamma_{\Delta}$. 
\vspace{4ex}

\noindent
{\bf 3. Universal abelian differential}
\vspace{2ex}

\noindent
{\it 3.1. Universal abelian differential.} 
Let $\Delta = (V, E, T)$ be a stable graph with 
$$
{\rm rank}_{\mathbb Z} H_{1}(\Delta, {\mathbb Z}) = g, 
$$
and the notation be as in 2.2. 
\vspace{2ex}

\noindent
{\bf Proposition 3.1.} 
\begin{it}
Let $\phi$ be a product $\phi_{h(1)} \cdots \phi_{h(l)}$ with 
$v_{-h(i)} = v_{h(i+1)}$ $(1 \leq i \leq l-1)$ 
which is reduced in the sense that $h(i) \neq -h(i+1)$ $(1 \leq i \leq l-1)$, 
and put $y_{\phi} = y_{h(1)} \cdots y_{h(l)}$. 

{\rm (1)} 
One has 
$\displaystyle \phi(z) - x_{h(1)} \in y_{h(1)} 
\left( R_{\Delta} \left[ z, \prod_{h \in \pm E} (z - x_{h})^{-1} \right] [[ y_{e} \ (e \in E) ]] \right)$. 

{\rm (2)} 
If $a \in A_{\Delta}$ satisfies $a - x_{-h(l)} \in A_{\Delta}^{\times}$, 
then $\phi(a) - x_{h(1)} \in I$. 
Furthermore, if $a' - x_{-h(l)} \in A_{\Delta}^{\times}$, 
then $\phi(a) - \phi(a') \in (a - a') y_{\phi} A_{\Delta}^{\times}$. 

{\rm (3)} 
One has 
$\displaystyle 
\frac{d \phi(z)}{dz} \in y_{\phi} 
\left( R_{\Delta} \left[ \prod_{h \in \pm E} (z - x_{h})^{-1} \right] [[ y_{e} \ (e \in E) ]] \right)$. 
\end{it} 
\vspace{2ex}

\noindent
{\it Proof.}
Since the assertion (2) is proved in \cite[Lemma 1.2]{I2}, 
we will prove (1) and (3). 
Put 
$$
\phi = \left( \begin{array}{cc} a_{\phi} & b_{\phi} \\ c_{\phi} & d_{\phi} \end{array} \right). 
$$
Since 
$$
\left( \begin{array}{cc} \alpha & - \alpha \beta \\ 1 & - \beta \end{array} \right)
\left( \begin{array}{cc} \gamma & - \gamma \delta \\ 1 & - \delta \end{array} \right) 
= (\gamma - \beta) 
\left( \begin{array}{cc} \alpha & - \alpha \beta \\ 1 & - \delta \end{array} \right), 
$$
$a_{\phi}$, $b_{\phi}$, $c_{\phi}$ and $d_{\phi}$ are elements of $A_{\Delta}$ 
whose constant terms are $x_{h(1)} t$, $- x_{h(1)} x_{-h(l)} t$, $t$ and $- x_{-h(l)} t$ 
respectively,   
where 
$$
t = \frac{\prod_{s = 2}^{l} \left( x_{h(s)} - x_{-h(s-1)} \right)}
{\prod_{s = 1}^{l} \left( x_{h(s)} - x_{-h(s)} \right)} \in A_{\Delta}^{\times}. 
$$
Then $c_{\phi} z + d_{\phi} = t (z - x_{-h(l)}) + \cdots$, 
and hence  
$$
\phi(z) - x_{h(1)} \in 
R_{\Delta} \left[ z, \prod_{h \in \pm E} (z - x_{h})^{-1} \right] [[ y_{e} \ (e \in E) ]]. 
$$
In order to prove (1), we may assume that $l = 1$, 
and then $\phi(z) = x_{h(1)}$ under $y_{h(1)} = 0$. 
Therefore, the assertion (1) holds.  
The assertion (3) follows from 
$$
\frac{d \phi(z)}{dz} = \frac{\det(\phi)}{(c_{\phi} z + d_{\phi})^{2}} 
= \frac{\prod_{s = 1}^{l} \det(\phi_{h(s)})}{(c_{\phi} z + d_{\phi})^{2}} 
= \frac{\prod_{s = 1}^{l} y_{h(s)}}{(c_{\phi} z + d_{\phi})^{2}}, 
$$
and the above calculation. 
\ $\square$ 
\vspace{2ex} 

For a stable graph $\Delta = (V, E, T)$, 
we define universal abelian differentials on a generalized Tate curve ${\mathcal C}_{\Delta}$. 
Let 
$\Gamma_{\Delta} = {\rm Im} \left( \pi_{1}(\Delta, v_{b}) \rightarrow 
PGL_{2}(B_{\Delta}) \right)$ 
be the universal Schottky group as above. 
Then it is shown in \cite[Lemma 1.3]{I2} that 
each $\gamma \in \Gamma_{\Delta} - \{ 1 \}$ has its attractive (resp. repulsive) 
fixed points $\alpha$ (resp. $\alpha'$) in ${\mathbb P}^{1}_{B_{\Delta}}$ 
and its multiplier $\beta \in \sum_{e \in E} A_{\Delta} \cdot y_{e}$ which satisfy 
$$
\frac{\gamma(z) - \alpha}{z - \alpha} = \beta \frac{\gamma(z) - \alpha'}{z - \alpha'}. 
$$ 
Fix a set $\{ \gamma_{1},..., \gamma_{g} \}$ of generators of $\Gamma_{\Delta}$, 
and for each $\gamma_{i}$, 
denote by $\alpha_{i}$ (resp. $\alpha_{-i}$) its attractive (resp. repulsive) fixed points, 
and by $\beta_{i}$ its multiplier.  
Then under the assumption that there is no element of $\pm E \cap {\mathcal E}_{\infty}$ 
with terminal vertex $v_{b}$, 
for each $1 \leq i \leq g$, 
we define the associated {\it universal abelian differential of the first kind} as 
$$
\omega_{i} = 
\sum_{\gamma \in \Gamma_{\Delta} / \left\langle \gamma_{i} \right\rangle} 
\left( \frac{1}{z - \gamma(\alpha_{i})} - \frac{1}{z - \gamma(\alpha_{-i})} \right) dz. 
$$
Assume that 
$$
\{ h \in \pm E \cap {\mathcal E}_{\infty} \ | \ v_{h} = v_{b} \} = \emptyset,  \ \ 
\{ t \in T \ | \ v_{t} = v_{b} \} \neq \emptyset. 
$$ 
Then for each $t \in T$ with $v_{t} = v_{b}$ and $k > 1$, 
we define the associated {\it universal abeian differential of the second kind} as 
$$
\omega_{t, k} = \left\{ \begin{array}{ll} 
{\displaystyle \sum_{\gamma \in \Gamma_{\Delta}} \frac{d \gamma(z)}{(\gamma(z) - x_{t})^{k}}} 
& (x_{t} \neq \infty), 
\\
\\
{\displaystyle \sum_{\gamma \in \Gamma_{\Delta}} \gamma(z)^{k-2} d \gamma(z)} 
& (x_{t} = \infty). 
\end{array} \right. 
$$
Furthermore, 
put $T_{\infty} = \{ t \in T \cap {\mathcal E}_{\infty} \ | \ v_{t} = v_{b} \}$ 
whose cardinality is $0$ or $1$, 
and take a maximal subtree ${\mathcal T}_{\Delta}$ of $\Delta$, 
and for each $t \in T$, 
take the unique path $\rho_{t} = h(1) \cdots h(l)$ in ${\mathcal T}_{\Delta}$ 
from $v_{t}$ to $v_{b}$, 
and put $\phi_{t} = \phi_{h(l)} \cdots \phi_{h(1)}$. 
Then for each $t_{1}, t_{2} \in T$ with $t_{1} \neq t_{2}$, 
we define the associated {\it universal abelian differential of the the third kind} as 
$$
\omega_{t_{1}, t_{2}} = \sum_{\gamma \in \Gamma_{\Delta}} 
\left( \frac{d \gamma(z)}{\gamma(z) - \phi_{t_{1}}(x_{t_{1}})} - 
\frac{d \gamma(z)}{\gamma(z) - \phi_{t_{2}}(x_{t_{2}})} \right), 
$$
where $\phi_{t_{i}}(x_{t_{i}}) = \infty$ if $t_{i} \in T_{\infty}$. 
\vspace{2ex}

\noindent
{\bf Theorem 3.2.} 
\begin{it} 

{\rm (1)} 
For each $1 \leq i \leq g$, 
$\omega_{i}$ is a regular differential on ${\mathcal C}_{\Delta} \otimes_{A_{\Delta}} B_{\Delta}$. 

{\rm (2)} 
For each $t \in T$ with $v_{t} = v_{b}$ and $k > 1$, 
$\omega_{t, k}$ is a meromorphic differential on ${\mathcal C}_{\Delta} \otimes_{A_{\Delta}} B_{\Delta}$ 
which has only pole (of order $k$) at the point $p_{t}$ corresponding to $t$. 

{\rm (3)}
For each $t_{1}, t_{2} \in T$ such that $t_{1} \neq t_{2}$, 
$\omega_{t_{1}, t_{2}}$ is a meromorphic differential on 
${\mathcal C}_{\Delta} \otimes_{A_{\Delta}} B_{\Delta}$ 
which has only (simple) poles at the points $p_{t_{1}}$ (resp. $p_{t_{2}}$) 
corresponding to $t_{1}$ (resp. $t_{2}$) with residue $1$ (resp. $-1$). 

{\rm (4)} 
Take $x_{h}$ $(h \in \pm E \cup T)$ and $y_{e}$ $(e \in E)$ 
be complex numbers as in 2.2 (P4). 
Then $\omega_{i}, \omega_{t, k}, \omega_{t_{1}, t_{2}}$ are relative abelian differentials 
on the family ${\mathcal R}_{\Delta}$ of Riemann surfaces. 

\end{it} 
\vspace{2ex}

\noindent
{\it Proof.} 
By Proposition 3.1, 
$\omega_{i}$ are differentials on ${\mathcal P}_{T_{\Delta}}$, 
and for any $\delta \in \Gamma_{\Delta}$, 

\begin{eqnarray*}
\omega_{i}(\delta(z)) & = & 
\sum_{\gamma \in \Gamma_{\Delta} / \left\langle \gamma_{i} \right\rangle} 
\left( \frac{1}{\delta(z) - \gamma(\alpha_{i})} - \frac{1}{\delta(z) - \gamma(\alpha_{-i})} \right) 
d \delta(z) 
\\ 
& = & 
\sum_{\gamma \in \Gamma_{\Delta} / \left\langle \gamma_{i} \right\rangle} 
\left( \frac{\gamma(\alpha_{i}) - \gamma(\alpha_{-i})}
{(\delta(z) - \gamma(\alpha_{i})) (\delta(z) - \gamma(\alpha_{-i}))} \right) 
d \delta(z) 
\\ 
& = & 
\sum_{\gamma \in \Gamma_{\Delta} / \left\langle \gamma_{i} \right\rangle} 
\left( \frac{(\delta^{-1} \gamma)(\alpha_{i}) - (\delta^{-1} \gamma)(\alpha_{-i})}
{(z - (\delta^{-1} \gamma)(\alpha_{i})) (z - (\delta^{-1} \gamma)(\alpha_{-i}))} \right) dz 
\\ 
& = & 
\sum_{\gamma \in \Gamma_{\Delta} / \left\langle \gamma_{i} \right\rangle} 
\left( \frac{1}{z - (\delta^{-1} \gamma)(\alpha_{i})} - 
\frac{1}{z - (\delta^{-1} \gamma)(\alpha_{-i})} \right) dz 
\\ 
& = & 
\omega_{i}(z). 
\end{eqnarray*}
Therefore, by the construction of ${\mathcal C}_{\Delta}$ reviewed in 2.2, 
$\omega_{i}$ give rise to differentials on ${\mathcal C}_{\Delta}$ 
which are regular outside $\bigcup_{e \in E} \{ y_{e} = 0 \}$, 
and hence the assertions (1) follows. 
One can prove the assertions (2), (3) similarly, 
and hence it is enough to prove (4). 
As is stated in 2.1, 
$R_{\Gamma}$ is given by the quotient space $\Omega_{\Gamma}/\Gamma$. 
Under the assumption on complex numbers $x_{h}$ and $y_{e}$, 
it is shown in \cite{Sc} that $\sum_{\gamma \in \Gamma} |\gamma'(z)|$ 
is uniformly convergent on any compact subset in 
$\Omega_{\Gamma} - \cup_{\gamma \in \Gamma} \gamma(\infty)$, 
and hence the assertion (4) holds for $\omega_{t, k}$, $\omega_{t_{1}, t_{2}}$. 
If $a \in \Omega_{\Gamma} - \cup_{\gamma \in \Gamma} \gamma(\infty)$, 
then $\lim_{n \rightarrow \infty} \gamma_{i}^{\pm n}(a) = \alpha_{\pm i}$, 
and hence 
\begin{eqnarray*}
d \left( \int_{a}^{\gamma_{i}(a)} \sum_{\gamma \in \Gamma} 
\frac{d \gamma(\zeta)}{\gamma(\zeta) - z} \right) 
& = & 
\sum_{\gamma \in \Gamma} 
\left( \frac{1}{z - (\gamma \gamma_{i})(a)} - \frac{1}{z - \gamma(a)} \right) dz 
\\ 
& = & 
\sum_{\gamma \in \Gamma / \langle \gamma_{i} \rangle} 
\lim_{n \rightarrow \infty} 
\left( \frac{1}{z - (\gamma \gamma_{i}^{n})(a)} - 
\frac{1}{z - (\gamma \gamma_{i}^{-n})(a)} \right) dz 
\\ 
& = & 
\omega_{i}(z). 
\end{eqnarray*}
Therefore, 
$\omega_{i}$ is absolutely and uniformly convergent 
on any compact subset in $\Omega_{\Gamma}$, 
and hence is an abelian differential on $\Omega_{\Gamma}/\Gamma$. 
\ $\square$ 
\vspace{2ex}

\noindent
{\it 3.2. Stability of universal abelian differentials.} 
For a vertex $v \in V$, 
denote by $C_{v}$ the corresponding irreducible component of 
${\mathcal C}_{\Delta} \otimes_{A_{\Delta}} R_{\Delta}$. 
Then $P_{v} = {\mathbb P}^{1}_{R_{\Delta}}$ is the normalization of $C_{v}$. 
\vspace{2ex}

\noindent
{\bf Theorem 3.3.} 
\begin{it} 

{\rm (1)} 
For each $1 \leq i \leq g$, 
let $\phi_{h_{i}(1)} \cdots \phi_{h_{i}(l_{i})}$ $(h_{i}(j) \in \pm E)$ be the unique reduced product 
such that $v_{-h_{i}(j)} = v_{h_{i}(j+1)}$ and $h_{i}(1) \neq - h_{i}(l_{i})$ which is conjugate to $\gamma_{i}$. 
Then for each $v \in V$, 
the pullback $\left( \omega_{i}|_{C_{v}} \right)^{*}$ of $\omega_{i}|_{C_{v}}$ to $P_{v}$ 
is given by 
$$
\left( \sum_{v_{h_{i}(j)} = v} \frac{1}{z - x_{h_{i}(j)}} - 
\sum_{v_{-h_{i}(k)} = v} \frac{1}{z - x_{-h_{i}(k)}} \right) dz. 
$$

{\rm (2)}
For each $v \in V$, 
$\left( \omega_{t, k}|_{C_{v}} \right)^{*}$ is given by 
$\displaystyle \frac{dz}{(z - x_{t})^{k}}$ if $v = v_{t}$, 
and is $0$ otherwise. 

{\rm (3)} 
Denote by $\rho_{t_{j}} = h_{j}(1) \cdots h_{j}(l_{j})$ 
the unique path from $v_{t_{j}}$ $(t_{j} \in T)$ to $v_{b}$ in ${\mathcal T}_{\Delta}$. 
Then for each $v \in V$, 
$\left( \omega_{t_{1}, t_{2}}|_{C_{v}} \right)^{*}$ is given by 
$$
\left( \sum_{v_{h} = v} \frac{1}{z - x_{h}} - 
\sum_{v_{-k} = v} \frac{1}{z - x_{-k}} \right) dz, 
$$
where $h, k$ runs through 
$\left\{ t_{1}, h_{1}(1),..., h_{1}(l_{1}), -t_{2}, -h_{2}(1),..., -h_{2}(l_{2}) \right\}$. 
\end{it} 
\vspace{2ex}

\noindent
{\it Proof.} 
For the proof of (1), 
we may assume that 
$\gamma_{i} = \phi_{h_{i}(1)} \cdots \phi_{h_{i}(l_{i})}$. 
Let $\gamma$ be an element of $\Gamma_{\Delta}$. 
Then by Proposition 3.1 (2), putting $y_{e} = 0$ $(e \in E)$, 
$$
\frac{1}{z - \gamma(\alpha_{i})} - \frac{1}{z - \gamma(\alpha_{-i})} = 
\frac{\gamma(\alpha_{i}) - \gamma(\alpha_{-i})}
{(z - \gamma(\alpha_{i}))(z - \gamma(\alpha_{-i}))}  
$$
becomes 
$$
\frac{1}{z - x_{h_{i}(j)}} - \frac{1}{z - x_{-h_{i}(j-1)}}
$$ 
if $j \in \{ 1,..., l_{i} \}$ $\left( \mbox{$h_{i}(j-1) := h_{i}(l_{i})$ when $j = 1$} \right)$, 
$v(h_{i}(j)) = v$, $\phi_{h_{i}(j)} \phi_{h_{i}(j+1)} \cdots \phi_{h_{i}(l_{i})}$ belongs to 
$\gamma \langle \gamma_{i} \rangle$, 
and becomes $0$ otherwise. 
Therefore, the assertion follows from the definition of $\omega_{i}$. 

The assertion (2) follows from Proposition 3.1 (1) and (3), 
and the assertion (3) can be shown in the same way as above. 
\ $\square$ 
\vspace{2ex}

A  (regular or meromorphic) global section of the dualizing sheaf on a stable curve 
is called a {\it stable differential} (cf. \cite{DM}). 
\vspace{2ex}

\noindent
{\bf Theorem 3.4.} 
\begin{it} 

{\rm (1)} 
For each $1 \leq i \leq g$, 
$\omega_{i}$ is a regular stable differential on ${\mathcal C}_{\Delta}/A_{\Delta}$, 
namely an element of the space 
$H^{0} \left( {\mathcal C}_{\Delta}, \omega_{{\mathcal C}_{\Delta}/A_{\Delta}} \right)$ 
of global sections of the dualizing sheaf $\omega_{{\mathcal C}_{\Delta}/A_{\Delta}}$ 
on ${\mathcal C}_{\Delta}/A_{\Delta}$. 
Furthermore, $\{ \omega_{i} \}_{1 \leq i \leq g}$ gives a basis of 
$H^{0} \left( {\mathcal C}_{\Delta}, \omega_{{\mathcal C}_{\Delta}/A_{\Delta}} \right)$. 

{\rm (2)} 
For each $t \in T$ with $v_{t} = v_{b}$ and $k > 1$, 
$\omega_{t, k}$ is a meromorphic stable differential on ${\mathcal C}_{\Delta}/A_{\Delta}$ 
which has only pole (of order $k$) at the point $p_{t}$ corresponding to $t$. 

{\rm (3)}
For each $t_{1}, t_{2} \in T$ with $t_{1} \neq t_{2}$, 
$\omega_{t_{1}, t_{2}}$ is a meromorphic stable differential on ${\mathcal C}_{\Delta}/A_{\Delta}$ 
which has only (simple) poles at the points $p_{t_{1}}$ (resp. $p_{t_{2}}$) 
corresponding to $t_{1}$ (resp. $t_{2}$) with residue $1$ (resp. $-1$). 

\end{it} 
\vspace{2ex}

\noindent
{\it Proof.} 
We only show that the latter assertion in (1) 
since the remains follow from Theorems 3.2 and 3.3. 
For each $i$, take the product $\phi_{h_{i}(1)} \cdots \phi_{h_{i}(l_{i})}$ of 
$\phi_{h}$ $(h \in \pm E)$ as in Theorem 3.3 (1). 
Then by this theorem, 
$\omega_{i}|_{C_{v_{h_{i}(1)}}}$ has simple pole at $x_{h_{i}(1)}$ with residue $1$, 
and hence $\{ \omega_{i} \}_{1 \leq i \leq g}$ gives a basis of 
$H^{0} \left( {\mathcal C}_{\Delta}, \omega_{{\mathcal C}_{\Delta}/A_{\Delta}} \right)$ 
since $\{ \gamma_{i} \}$ is a basis of 
$H_{1}(\Delta, {\mathbb Z}) \cong 
\Gamma_{\Delta}/\left[ \Gamma_{\Delta}, \Gamma_{\Delta} \right]$. 
\ $\square$ 
\vspace{2ex}

\noindent
{\it Remark 3.5.} 
For a subset $E'$ of the set $E$ of edges in the stable graph $\Delta$, 
let $\Delta' = (V', E', T')$ be the associated graph obtained from $\Delta$, 
and denote by ${\mathcal C}_{\Delta'}$ the stable curve given by putting $y_{e} = 0$ $(e \in E')$. 
Let ${\mathcal R}_{\Delta}$ be a family of Schottky uniformized Riemann surfaces as in 2.2 (P4), 
and let $c_{e}$ be a cycle corresponding to $e \in E - E'$ which is oriented by 
the right-hand rule for $\gamma_{i}$. 
Then by Theorem 3.3 (1), 
the restriction of $\omega_{i}$ to each irreducible component of ${\mathcal R}_{\Delta'}$ 
is characterized analytically as follows: 
\begin{itemize}

\item 
its integrals along $c_{e}$ $(e \in E - E')$ are $2 \pi \sqrt{-1}$ 
if $e = |h_{i}(j)|$ for some integer $1 \leq j \leq l_{i}$, and are $0$ otherwise, 
and its residues at $x_{h_{i}(j)}$ (resp. $x_{h_{i}(k)}$) are $1$ (resp. $-1$) 
for some integers $1 \leq j, k \leq l_{i}$, and are $0$ otherwise. 
\end{itemize} 
By results in Section 5 below, 
this characterization holds for any deformation of a stable complex curve 
whose irreducible components are not necessarily projective lines. 
This modification gives a more explicit formula than \cite[Corollary 4.6]{HN}, 
and one has similar statements on $\omega_{t,k}$ and $\omega_{t_{1}, t_{2}}$. 
\vspace{2ex}

\noindent
{\it 3.3. Period.} 
\vspace{2ex}

\noindent
{\bf Theorem 3.6.} 
\begin{it} 

{\rm (1)} 
The Jacobian of ${\mathcal C}_{\Delta}$ becomes the Mumford abelian scheme 
\cite{Mu2} whose multiplicative periods is given in \cite[3.10]{I2} 
as $P_{ij} = \prod_{\gamma} \psi_{ij}(\gamma)$ $(1 \leq i, j \leq g)$, 
where $\gamma$ runs through all representatives of 
$\langle \gamma_{i} \rangle \backslash \Gamma_{\Delta} / \langle \gamma_{j} \rangle$ 
and 
$$
\psi_{ij}(\gamma) = \left\{ \begin{array}{ll} 
\beta_{i} & (\mbox{$i = j$ and $\gamma \in \langle \gamma_{i} \rangle$}), 
\\
{\displaystyle \frac{(\alpha_{i} - \gamma(\alpha_{j}))(\alpha_{-i} - \gamma(\alpha_{-j}))}
{(\alpha_{-i} - \gamma(\alpha_{j}))(\alpha_{i} - \gamma(\alpha_{-j}))}} 
& (\mbox{otherwise}). 
\end{array} \right. 
$$ 

{\rm (2)} 
Let $R$ be a Riemann surface obtained from ${\mathcal C}_{\Delta}$ as in 2.2 (P4), 
and take counterclockwise oriented small closed paths $a_{i}$ $(1 \leq i \leq g)$ in $R$ 
surrounding $\alpha_{i}$, 
and closed paths $b_{i}$ $(1 \leq i \leq g)$ in $R$ corresponding to $\gamma_{i}$. 
Denote by the same symbols the abelian differentials of the first kind and 
the multiplicative periods for $R$ given by the above $\omega_{j}$ $(1 \leq j \leq g)$ and $P_{ij}$ $(1 \leq i, j \leq g)$ respectively. 
Then 
$$
\int_{a_{i}} \omega_{j} = 2 \pi \sqrt{-1} \delta_{ij}, \ \ 
\exp \left( \int_{b_{i}} \omega_{j} \right) = P_{ij}, 
$$  
where $\delta_{ij}$ denotes the Kronecker delta.
\end{it}
\vspace{2ex} 

\noindent
{\it Proof.} 
Since the assertion (1) is shown in \cite[Theorem 3.13]{I2} 
and the first formula of the assertion (2) follows from Remark 3.5, 
we will prove the second formula of (2) (cf. \cite{Sc} and \cite[Theorem 2]{MD}). 
Let $\Gamma$ be a Schottky group over ${\mathbb C}$ obtained from 
$\Gamma_{\Delta}$ as in 2.2 (P4). 
Then as is seen in 2.1, 
$R = R_{\Gamma}$ is obtained as the quotient space of 
the discontinuity set $\Omega_{\Gamma}$ under the action of $\Gamma$. 
Take a point $a$ on $\Omega_{\Gamma}$. 
Then under the condition that $i = j$ and $\gamma \in \langle \gamma_{i} \rangle$,  
$$
\frac{\left( \gamma_{i}(a) - \gamma(\alpha_{i}) \right)
\left( a - \gamma(\alpha_{-i}) \right)}
{\left( a - \gamma(\alpha_{i}) \right)
\left( \gamma_{i}(a) - \gamma(\alpha_{-i}) \right)} 
= \beta_{i} = \psi_{ii}(\gamma). 
$$
Since $\lim_{n \rightarrow \infty} \gamma_{i}^{\pm n}(a) = \alpha_{\pm i}$, 
if $i \neq j$ or $\gamma \not\in \langle \gamma_{i} \rangle$, then 
\begin{eqnarray*}
\lefteqn{
\prod_{n \in {\mathbb Z}} 
\left( \frac{\left( \gamma_{i}(a) - \gamma_{i}^{-n} \gamma(\alpha_{j}) \right)
\left( a - \gamma_{i}^{-n} \gamma(\alpha_{-j}) \right)}
{\left( a - \gamma_{i}^{-n} \gamma(\alpha_{j}) \right)
\left( \gamma_{i}(a) - \gamma_{i}^{-n} \gamma(\alpha_{-j}) \right)} \right) 
}
\\
& = & 
\prod_{n \in {\mathbb Z}} 
\left( \frac
{\left( \gamma_{i}^{n+1}(a) - \gamma(\alpha_{j}) \right)
\left( \gamma_{i}^{n}(a) - \gamma(\alpha_{-j}) \right)}
{\left( \gamma_{i}^{n}(a) - \gamma(\alpha_{j}) \right)
\left( \gamma_{i}^{n+1}(a) - \gamma(\alpha_{-j}) \right)} \right) 
\\ 
& = & 
\frac{(\alpha_{i} - \gamma(\alpha_{j}))(\alpha_{-i} - \gamma(\alpha_{-j}))}
{(\alpha_{-i} - \gamma(\alpha_{j}))(\alpha_{i} - \gamma(\alpha_{-j}))} 
\\ 
& = & 
\psi_{ij}(\gamma). 
\end{eqnarray*}
Therefore, 
$$
\exp \left( \int_{b_{i}} \omega_{j} \right) 
= 
\prod_{\gamma \in \langle \gamma_{i} \rangle \backslash \Gamma/ \langle \gamma_{j} \rangle} 
\psi_{ij}(\gamma) 
= P_{ij} 
$$
which completes the proof. 
\ $\square$ 
\vspace{4ex}

\noindent
{\bf 4. Universal Mumford curve, differentials and periods} 
\vspace{2ex}

\noindent 
{\it 4.1. Universal Mumford curve.} 
Let $\Delta = (V, E, T)$ be a stable graph which is not trivalent. 
Then there exists a vertex $v_{0} \in V$ which has at least $4$ branches. 
Take two elements $h_{1}, h_{2}$ of ${\mathcal E}$ such that 
$h_{1} \neq h_{2}$ and $v_{h_{1}} = v_{h_{2}} = v_{0}$, 
and let $\Delta' = (V', E', T')$ be a stable graph obtained from $\Delta$ 
by replacing $v_{0}$ with an oriented (nonloop) edge $h_{0}$ such that 
$v_{h_{1}} = v_{h_{2}} = v_{h_{0}}$ and that $v_{h} = v_{-h_{0}}$ 
for any $h \in \pm E \cup T - \{ h_{1}, h_{2} \}$ with $v_{h} = v_{0}$. 
Put $e_{i} = |h_{i}|$ for $i = 0, 1, 2$. 
Then we have the following identifications: 
$$
V = V' - \{ v_{-h_{0}} \} \ (\mbox{in which $v_{0} = v_{h_{0}}$}), \ E = E' - \{ e_{0} \}, \ T = T'. 
$$ 

\noindent 
{\bf Theorem 4.1} 
\begin{it} 

{\rm (1)} 
The generalized Tate curves ${\mathcal C}_{\Delta}$ and ${\mathcal C}_{\Delta'}$ 
associated with $\Delta$ and $\Delta'$ respectively are isomorphic 
over $R_{\Delta'} [[ s_{e} \ (e \in E' - \{ e_{0} \}) ]] [s_{e_{0}}^{-1}]$, 
where 
$$
\frac{x_{h_{1}} - x_{h_{2}}}{s_{e_{0}}}, \ \ 
\frac{y_{e_{i}}}{s_{e_{0}} s_{e_{i}}} \ (\mbox{$i = 1, 2$ with $h_{i} \not\in T$}), \ \
\frac{y_{e}}{s_{e}} \ (e \in E - \{ e_{1}, e_{2} \}) 
$$
belong to $(A_{\Delta'})^{\times}$ if $h_{1} \neq - h_{2}$, 
and 
$$
\frac{x_{h_{1}} - x_{h_{2}}}{s_{e_{0}}}, \ \ 
\frac{y_{e}}{s_{e}} \ (e \in E) 
$$
belong to $(A_{\Delta'})^{\times}$ if $h_{1} = - h_{2}$. 

{\rm (2)} 
The assertion (1) holds in the category of complex geometry 
when $x_{h_{1}} - x_{h_{2}}, y_{e}$ and $s_{e}$ are taken to be 
sufficiently small complex numbers with $x_{h_{1}} \neq x_{h_{2}}, s_{e_{0}} \neq 0$. 

\end{it}
\vspace{2ex}

\noindent 
{\it Proof.} 
First, we prove the assertion (1). 
Denote by $t_{h}$ the moduli parameters of ${\mathcal C}_{\Delta'}$ corresponding to 
$h \in \pm E' \cup T'$.  
Then by \cite[Lemma 1.2]{I1}, 
$\phi_{-h_{0}}(t_{h_{1}}) - \phi_{-h_{0}}(t_{h_{2}})$ belongs to 
$s_{e_{0}} \cdot (A_{\Delta'})^{\times}$, 
and hence 
$$
{\mathcal C}_{\Delta'} \otimes_{A_{\Delta'}} 
\left( R_{\Delta'} [[ s_{e} ]] [s_{e_{0}}^{-1}] \right)
$$ 
gives a universal deformation of a universal degenerate curve with dual graph $\Delta$. 
Then by the universality of generalized Tate curves, 
there exists an injective homomorphism 
$A_{\Delta} \hookrightarrow R_{\Delta'} [[ s_{e} ]] [s_{e_{0}}^{-1}]$
which gives rise to an isomorphism ${\mathcal C}_{\Delta} \cong {\mathcal C}_{\Delta'}$. 
Under this homomorphism, 
\begin{eqnarray*}
\lefteqn{
\left( P_{v_{-h_{0}}}; \phi_{-h_{0}}(t_{h_{1}}), \phi_{-h_{0}}(t_{h_{2}}), 
t_{h} \ (v_{h} = v_{-h_{0}}, h \neq h_{0}) \right) 
} 
\\ 
& \cong & 
\left( P_{v_{0}}; x_{h_{1}}, x_{h_{2}}, x_{h} \ (v_{h} = v_{0}, h \neq h_{1}, h_{2}) \right), 
\end{eqnarray*} 
and hence $x_{h_{1}} - x_{h_{2}} \in s_{e_{0}} \cdot (A_{\Delta'})^{\times}$. 
Furthermore, when $h_{1} \neq -h_{2}$, 
the deformation parameters of 
$$
\left( P_{v_{-h_{0}}}; \phi_{-h_{0}}(t_{h_{1}}), \phi_{-h_{0}}(t_{h_{2}}), 
t_{h} \ (v_{h} = v_{-h_{0}}, h \neq h_{0}) \right) 
$$
corresponding to $h_{i} \cdot (-h_{0})$ $(i = 1, 2)$ are $y_{h_{i}}$, 
and hence by \cite[Proposition 1.3]{I1}, 
$y_{h_{i}} \in \left( s_{h_{0}} \cdot s_{h_{i}} \right) \cdot  (A_{\Delta'})^{\times}$. 
When $h_{1} = - h_{2}$, 
the deformation parameters of 
$$
\left( P_{v_{-h_{0}}}; \phi_{-h_{0}}(t_{h_{1}}), \phi_{-h_{0}}(t_{h_{2}}), 
t_{h} \ (v_{h} = v_{-h_{0}}, h \neq h_{0}) \right) 
$$
corresponding to $h_{0} \cdot h_{1} \cdot (-h_{0})$ is $y_{h_{1}}$, 
and hence $y_{h_{1}} \in s_{h_{1}} \cdot  (A_{\Delta'})^{\times}$. 

The assertion (2) follows from 2.2 (P4). 
\ $\square$ 
\vspace{2ex}

\noindent
{\it Remark 4.2} 
Let $\Gamma_{\Delta}$ and $\Gamma_{\Delta'}$ denote the universal Schottky groups 
associated with $\Delta$ and $\Delta'$ respectively. 
Since ${\mathcal C}_{\Delta} \cong {\mathcal C}_{\Delta'}$, 
by \cite[Corollary 4.11]{Mu1}, there exists a conjugation isomorphism 
$\varphi : \Gamma_{\Delta} \stackrel{\sim}{\rightarrow} \Gamma_{\Delta'}$, 
and hence the multiplier of an element of $\Gamma_{\Delta}$ 
and the cross-ratio of the fixed points of four elements of $\Gamma_{\Delta}$ 
are invariant under $\varphi$. 
Therefore, one can calculate the precise formula connecting 
the moduli and deformation parameters for ${\mathcal C}_{\Delta}$ and ${\mathcal C}_{\Delta'}$. 
For the detail, see the proof of \cite[Theorem 1]{I2}. 
\vspace{2ex}   
 
For nonnegative integers $g, n$ such that $2g - 2 + n > 0$, 
denote by $\overline{\mathcal M}_{g,n}$ the moduli stack over ${\mathbb Z}$ 
of stable $n$-marked curves of genus $g$ (cf. \cite{DM, Kn, KnM}). 
Then by definition, there exist the universal stable marked curve ${\mathcal C}_{g,n}$ 
over $\overline{\mathcal M}_{g,n}$, 
and the associated curve ${\mathcal C}_{g}$ 
obtained by forgetting marked points on ${\mathcal C}_{g,n}$. 
\vspace{2ex}

\noindent 
{\bf Theorem 4.3.} 
\begin{it} 
There exists a deformation space ${\mathcal A}_{g, n}$ of 
all $n$-marked degenerate curves of genus $g$, 
and an $n$-marked stable curve of genus $g$ over ${\mathcal A}_{g, n}$ 
whose fiber by the canonical morphism ${\rm Spec}(A_{\Delta}) \rightarrow {\mathcal A}_{g, n}$ 
becomes the generalized Tate curve ${\mathcal C}_{\Delta}$ 
for each stable graph $\Delta$ of $(g, n)$-type. 
\end{it} 
\vspace{2ex}

\noindent 
{\it Proof.} 
Let $\Delta = (V, E, T)$ be a stable graph of $(g, n)$-type, 
and take a system of coordinates on $P_{v} = {\mathbb P}^{1}_{R_{\Delta}}$ $(v \in V)$ 
such that $x_{h} = \infty$ $(h \in {\mathcal E}_{\infty})$ and that 
$\{ 0, 1 \} \subset P_{v}$ is contained in the set of points 
given by $x_{h}$ $(h \in {\mathcal E}$ with $v_{h} = v)$. 
Under this system of coordinates, 
one has the generalized Tate curve ${\mathcal C}_{\Delta}$ whose closed fiber 
${\mathcal C}_{\Delta} \otimes_{A_{\Delta}} R_{\Delta}$ gives a family of degenerate curves 
over the open subspace of 
$$
S_{\Delta} = \left\{ \left. \left( p_{h} \in P_{v_{h}} \right)_{h \in \pm E \cup T} \ \right| \ 
p_{h} \neq p_{h'} \ (h \neq h', v_{h} = v_{h'}) \right\} 
$$ 
defined as $p_{e} \neq p_{-e}$ for nonloop edges $e \in E$. 
Therefore, taking another system of coordinates on $P_{v}$ obtained by mutual changes 
of $0, 1, \infty$ and comparing the associated generalized Tate curves 
with the original ${\mathcal C}_{\Delta}$ as in Theorem 4.1, 
${\mathcal C}_{\Delta}$ can be extended over the deformation space of 
all marked degenerate curves with dual graph $\Delta$. 
Since two stable graphs of $(g, n)$-type can be translated by a combination of 
replacements $\Delta \leftrightarrow \Delta'$ given in 3.1, 
one can define a scheme ${\mathcal A}_{g, n}$ obtained by gluing 
${\rm Spec}(A_{\Delta})$ ($\Delta$: stable graphs of $(g, n)$-type) 
along the isomorphism given in Theorem 4.1. 
Then ${\mathcal A}_{g, n}$ is regarded as the deformation space of 
all $n$-marked degenerate curves of genus $g$ over which 
there exists an $n$-marked stable curve of genus $g$ obtained by gluing ${\mathcal C}_{\Delta}$. 
\ $\square$  
\vspace{2ex} 

\noindent 
{\it Definition 4.4.} 
We call the above $n$-marked stable curve of genus $g$ over ${\mathcal A}_{g, n}$  
the $n$-marked {\it universal Mumford curve} of genus $g$ 
which is the fiber of ${\mathcal C}_{g, n}$ by the canonical morphism 
${\mathcal A}_{g, n} \rightarrow \overline{\mathcal M}_{g, n}$. 
By 2.2 (P2) and (P4), 
one can see that this universal Mumford curve gives rise to 
all $n$-marked Mumford curves of genus $g$, 
and to all $n$-marked Riemann surfaces of genus $g$ close to degenerate curves. 
\vspace{2ex}

\noindent 
{\it Remark 4.5.} 
Gerritzen-Herrlich \cite{GH} introduced the extended Schottky space 
$\overline{S}_{g}$ of genus $g > 1$ 
as the fine moduli space of stable complex curves of genus $g$ with Schottky structure. 
For integers $g, n$ as above, 
one can consider the extended Schottky space $\overline{S}_{g,n}$ 
for stable $n$-marked complex curves of genus $g$ with Schottky structure. 
Then by the result of Koebe referred in 2.1, 
$\overline{S}_{g, n}/{\rm Out}(F_{g})$ becomes a covering space of 
the moduli space of $n$-marked stable complex curves of genus $g$ 
which is also the complex analytic space $\overline{\mathcal M}_{g, n}^{\rm an}$ 
associated with $\overline{\mathcal M}_{g, n}$. 
Furthermore, the $n$-marked universal Mumford curve of genus $g$ can be 
analytically  extended to the universal family of marked stable complex curves 
over $\overline{S}_{g,n}/{\rm Out}(F_{g})$. 
\vspace{2ex}

\noindent
{\it 4.2. Differentials and periods of the universal Mumford curve.}   
Denote by $\Delta_{0} = (V_{0}, E_{0}, T_{0})$ the stable graph of $(g, n)$-type 
consisting of one vertex and $g$ loops, 
and fix generators $\rho_{1},..., \rho_{g}$ of $\pi_{1}(\Delta_{0})$. 
For each stable graph $\Delta = (V, E, T)$ of $(g, n)$-type, 
$\Delta_{0}$ is obtained from $\Delta$ by contracting some nonloop edges in $E$. 
Then there exists a unique ${\mathbb Z}$-basis $\{ b_{i} \}_{1 \leq i \leq g}$ 
of $H_{1}(\Delta, {\mathbb Z})$ corresponding to 
$$
\left\{ [\rho_{i}] \right\}_{1 \leq i \leq g} \in H_{1}(\Delta_{0}, {\mathbb Z}), 
$$
and $T$ is identified with $T_{0}$.  
Therefore, by Theorem 3.4, 
there exist associated stable differentials 
$$
\omega_{i} \ (i = 1,..., g), \ \ \omega_{t, k} \ (t \in T_{0}, k > 1), \ \ 
\omega_{t_{1}, t_{2}} \ (t_{i} \in T_{0}, t_{1} \neq t_{2}) 
$$
on ${\mathcal C}_{\Delta}$ of the first, second, third kind respectively. 
Furthermore, $b_{i}$ $(i = 1,..., g)$ give rise to homology cycles 
on members of ${\mathcal R}_{\Delta}$ given in 2.2 (P4), 
and on the universal family of stable complex curves over 
the extended Schottky space $\overline{S}_{g,n}$ by the analytic extension. 
We also denote by $b_{i}$ these homology cycles. 
\vspace{2ex}

\noindent 
{\bf Theorem 4.6.} 
\begin{it} 

{\rm (1)} 
The differentials $\omega_{i}$, $\omega_{t, k}$, $\omega_{t_{1}, t_{2}}$ on ${\mathcal C}_{\Delta}$ 
are glued to stable differentials on ${\mathcal C}_{g}/{\mathcal A}_{g,n}$ 
which we call the universal differentials of the first kind, second kind, third kind, 
and denote by 
$\overline{\omega}_{i}$, $\overline{\omega}_{t, k}$, $\overline{\omega}_{t_{1}, t_{2}}$ respectively. 

{\rm (2)} 
The universal differentials $\overline{\omega}_{i}$ $(i = 1,..., g)$ on 
${\mathcal C}_{g}/{\mathcal A}_{g,n}$ of the first kind make a basis of the sheaf 
${\mathcal H}^{0} \left( \omega_{{\mathcal C}_{g}/{\mathcal A}_{g, n}} \right)$ 
which consists of sections of the relative stable differentials 
on ${\mathcal C}_{g}/{\mathcal A}_{g,n}$. 
Furthermore, $\overline{\omega}_{i}$ are analytically extended to 
stable differentials on the universal family of stable complex curves over $\overline{S}_{g,n}$ 
which we denote by the same symbols. 
\end{it}
\vspace{2ex} 

\noindent 
{\it Proof.} 
First, we prove (1). 
As is stated in Remark 4.2, 
the isomorphism ${\mathcal C}_{\Delta} \cong {\mathcal C}_{\Delta_{0}}$ over $B_{\Delta}$ 
considered in Theorem 3.1 corresponds uniquely to an isomorphism 
$\Gamma_{\Delta} \cong \Gamma_{\Delta_{0}}$. 
Therefore, under this isomorphism ${\mathcal C}_{\Delta} \cong {\mathcal C}_{\Delta_{0}}$, 
the differentials $\omega_{i}$, $\omega_{t, k}$, $\omega_{t_{1}, t_{2}}$ on ${\mathcal C}_{\Delta}$ 
are mapped to those on ${\mathcal C}_{\Delta_{0}}$, 
and hence can be glued to differentials on ${\mathcal C}_{g}/{\mathcal A}_{g,n}$.   

Second, we prove the latter assertion of (2) 
since the former one follows from Theorem 4.6 (1). 
For a stable $n$-marked complex curve $C$ of genus $g$ with Schottky structure, 
as is stated in \cite{GH}, 
one can take a cut system $\{ a_{1},..., a_{g} \}$ consisting of disjoint oriented simple loops 
in $C$ such that the intersection numbers $a_{i} \cdot b_{j}$ 
are the Kronecker delta $\delta_{ij}$. 
Then there exist uniquely regular stable differentials $\omega_{C, i}$ $(i = 1,..., g)$ on $C$ 
such that 
$$
\displaystyle \int_{a_{j}} \omega_{C, i} = 2 \pi \sqrt{-1} \delta_{ij}. 
$$
Therefore, moving $C$ on $\overline{S}_{g, n}$, 
$\omega_{C, i}$ form the analytic extension of $\overline{\omega}_{i}$.  
\ $\square$
\vspace{2ex}

\noindent 
Denote by ${\mathcal B}_{g, n}$ the maximal open subscheme of  ${\mathcal A}_{g, n}$ 
over which the associated marked curves are smooth, 
and by $S_{g, n}$ the (ordinary) Schottky space 
which is defined as the open subspace of $\overline{S}_{g, n}$ 
classifying $n$-marked Riemann surfaces of genus $g$ with Schottky structure. 
\vspace{2ex}

\noindent 
{\bf Theorem 4.7.} 
\begin{it} 
There exist ${\mathcal P}_{ij} \in {\mathcal O}^{\times}_{{\mathcal B}_{g, n}}$ 
$(1 \leq i, j \leq g)$ which give the multiplicative periods of Mumford curves 
over $p$-adic fields obtained from ${\mathcal C}_{g}/{\mathcal B}_{g, n}$ as in 2.2 (P2). 
Furthermore, ${\mathcal P}_{ij}$ can be analytically continued to regular functions 
on $S_{g, n}$ which become $\displaystyle \exp \left( \int_{b_{i}} \overline{\omega}_{j} \right)$, 
where $\overline{\omega}_{j}$ are given in Theorem 4.6. 
We call ${\mathcal P}_{ij}$ the universal periods. 
\end{it}
\vspace{2ex}

\noindent 
{\it Proof.} 
For a stable graph $\Delta$ of $(g, n)$-type, 
take generators $\gamma_{1},..., \gamma_{g}$ of $\Gamma_{\Delta}$ such that 
$[\gamma_{i}] = b_{i}$ $(i = 1,..., g)$, 
and for each $\gamma_{i}$, 
denote by $\alpha_{i}$ (resp. $\alpha_{-i}$) its attractive (resp. repulsive) fixed points 
and by $\beta_{i}$ its multiplier. 
Then by \cite[Theorem 3.13]{I1}, 
the multiplicative periods $P_{ij}$ $(1 \leq i, j \leq g)$ of ${\mathcal C}_{\Delta}$ are 
elements of $B_{\Delta}^{\times}$ defined as $P_{ij} = \prod_{\gamma} \psi_{ij}(\gamma)$, 
where $\gamma$ runs through all representatives of 
$\langle \gamma_{i} \rangle \backslash \Gamma_{\Delta} / \langle \gamma_{j} \rangle$ 
and 
$$
\psi_{ij}(\gamma) = \left\{ \begin{array}{ll} 
\beta_{i} & (i = j, \ \gamma \in \langle \gamma_{i} \rangle), 
\\
{\displaystyle \frac{(\alpha_{i} - \gamma(\alpha_{j})) (\alpha_{-i} - \gamma(\alpha_{-j}))}
{(\alpha_{-i} - \gamma(\alpha_{j})) (\alpha_{i} - \gamma(\alpha_{-j}))}} 
& (\mbox{otherwise}). 
\end{array} \right. 
$$ 
This implies that $P_{ij}$ depend only on $\{ b_{i} \}_{1 \leq i \leq g}$, 
and they give rise to the multiplicative periods 
$\langle [\gamma_{i}], [\gamma_{j}] \rangle$ of Mumford curves 
over $p$-adic fields by \cite[Theorem 2]{MD}, 
and to $\displaystyle \exp \left( \int_{b_{i}} \omega_{j} \right)$ 
by the definition of $\omega_{j}$ (see also \cite[Theorem 7.5]{M}). 
Therefore, as in the proof of Theorem 5.1, 
$P_{ij}$ are seen to be glued to regular functions ${\mathcal P}_{ij}$ on ${\mathcal B}_{g, n}$ 
satisfying the required properties. 
\ $\square$ 
\vspace{4ex}

\noindent
{\bf 5. Degeneration of abelian differentials}  
\vspace{2ex}

\noindent
{\it 5.1. Irreducible degeneration.} 
We apply results in Section 3 to the study of abelian differentials 
on a general family of Riemann surfaces which stably degenerates to 
an irreducible singular complex curve. 

Let ${\mathcal R}$ be a family of Riemann surfaces of genus $g > 0$ 
with symplectic bases $\{ a_{i}, b_{i} \}_{1 \leq i \leq g}$ 
(of the first Betti homology groups of members in ${\mathcal R}$) 
and local coordinates $u$ at marked points $p$ 
which is defined over a small disk around $y = 0$ in ${\mathbb C}$. 
Here $u$ are centered at $p$ that is $u(p) = 0$, 
and $\{ a_{i}, b_{i} \}$ are taken as homotopy classes of closed paths in members in ${\mathcal R}$  
such that $(a_{i}, b_{j}) = \delta_{ij}, (a_{i}, a_{j}) = (b_{i}, b_{j}) = 0$. 
Assume that letting $y \rightarrow 0$ corresponds to pinching $a_{g}$ to a point, 
and then ${\mathcal R}$ degenerates to a singular complex curve obtained from 
a Riemann surface $R'$ of genus $g-1$ by identifying two points $p_{1}, p_{2}$ on $R'$, 
where $b_{g}$ becomes (up to homotopy) to a path from $p_{1}$ to $p_{2}$ 
in $R'/(p_{1} = p_{2})$. 
Denote by $\left( \{ a'_{i}, b'_{i} \}_{1 \leq i \leq g-1}, p, u \right)$ 
a symplectic basis and a local coordinate at a marked point of $R'$ obtained from 
$\left( \{ a_{i}, b_{i} \}_{1 \leq i \leq g-1}, p, u \right)$ by putting $y = 0$. 

Let $\{ \omega_{1},..., \omega_{g} \}$ be a unique basis of the space of relative abelian differentials 
of the first kind on ${\mathcal R}$ which is normalized 
for $\{ a_{i} \}_{1 \leq i \leq g}$ in the sense that 
$$
\int_{a_{i}} \omega_{j} = 2 \pi \sqrt{-1} \delta_{ij} \ \ (1 \leq i, j \leq g). 
$$
By the theorem of Riemann-Roch, for each positive integer $n$, 
there exists a unique relative abelian differential $\omega^{(n)}$ of the second kind 
with order $n+1$ on ${\mathcal R}$ which is normalized for $(\{ a_{i} \}_{1 \leq i \leq g}, p, u)$ 
in the sense that 
\begin{itemize}

\item 
$\omega^{(n)}$ is holomorphic outside $p$, 

\item 
$\displaystyle \omega^{(n)} - \frac{1}{u^{n+1}} du$ is holomorphic at $p$, 

\item 
$\displaystyle \int_{a_{i}} \omega^{(n)} = 0$ for any $i = 1,..., g$. 

\end{itemize} 
\vspace{1ex}

\noindent
{\bf Theorem 5.1.} 
\begin{it} 

{\rm (1)} 
Under $y \rightarrow 0$, 
$\{ \omega_{1},..., \omega_{g-1} \}$ becomes to the basis of the space of abelian differentials of the first kind 
on $R'$ which is normalized for $\{ a'_{i} \}_{1 \leq i \leq g-1}$, 
and $\omega_{g}$ becomes to the unique abelian differential $\omega'$ of the third kind on $R'$ 
which is normalized for $(\{ a'_{i} \}_{1 \leq i \leq g-1}, p_{1}, p_{2})$ in the sense that 
$\omega'$ has only simple poles at $p_{1}, p_{2}$ with residue $1, -1$ respectively and satisfies 
$$
\int_{a'_{i}} \omega' = 0 \ \ (1 \leq i \leq g-1). 
$$   

{\rm (2)} 
Under $y \rightarrow 0$, 
$\omega^{(n)}$ becomes to the unique abelian differential of the second kind with order $n+1$ 
on $R'$ normalized for $(\{ a'_{i} \}_{1 \leq i \leq g-1}, p, u)$. 
\end{it}
\vspace{2ex}

\noindent
{\it Proof.}
First, we prove the assertions (1) and (2) in the case when ${\mathcal R}$ is obtained from 
a generalized Tate curve as in 2.2 (P4). 
Let $\Delta = (V, E, T)$ be a stable graph such that 
${\rm rank}_{\mathbb Z} H_{1}(\Delta, {\mathbb Z}) = g$, $\sharp V = \sharp T = 1$ 
and that there is an oriented loop in $\pm E$ which we denote by $l$, 
where $\Delta$ is obtained from a stable graph $\Delta'$ with three tails $t, t_{1}, t_{2}$ 
by connecting $t_{1}, t_{2}$ with $l$. 
There exist generators $\gamma_{1},..., \gamma_{g}$ of $\Gamma_{\Delta}$ 
such that $\gamma_{1},..., \gamma_{g-1}$ are generators of $\Gamma_{\Delta'}$ 
and that $\gamma_{g}$ corresponds to $l$. 
We take a symplectic basis $\{ a_{i}, b_{i} \}_{1 \leq i \leq g}$ of 
a family ${\mathcal R}_{\Delta}$ of Riemann surfaces obtained from ${\mathcal C}_{\Delta}$ 
such that $b_{i}$ corresponds to $\gamma_{i}$ as in 2.2 (P4). 
Then by Theorem 3.4, the universal abelian differentials of the first kind 
$$
\omega_{i} = 
\sum_{\gamma \in \Gamma_{\Delta} / \left\langle \gamma_{i} \right\rangle} 
\left( \frac{1}{z - \gamma(\alpha_{i})} - \frac{1}{z - \gamma(\alpha_{-i})} \right) dz \ \ 
(1 \leq i \leq g), 
$$
on ${\mathcal C}_{\Delta}$ give the basis of the space of relative abelian differentials of the first kind 
on ${\mathcal R}_{\Delta}$ normalized for $\{ a_{i} \}_{1 \leq i \leq g}$. 
By Proposition 3.1 (2), 
$$
\left. \left( \frac{1}{z - \gamma(\alpha_{i})} - \frac{1}{z - \gamma(\alpha_{-i})} \right) 
\right|_{y_{g}=0} = 
\left. \frac{\gamma(\alpha_{i}) - \gamma(\alpha_{-i})}
{(z - \gamma(\alpha_{i}))(z - \gamma(\alpha_{-i}))} 
\right|_{y_{g}=0} = 0 
$$
if $i \neq g$, $\gamma \in \Gamma_{\Delta} - \Gamma_{\Delta'}$ or 
$i = g$, $\gamma$ has a reduced product $\gamma' \gamma_{g}^{n}$, 
where $\gamma' \in \Gamma_{\Delta} - \Gamma_{\Delta'}$ and $n \in {\mathbb Z}$. 
Therefore, 
$$
\omega_{i}|_{y_{g} = 0} = 
\sum_{\gamma \in \Gamma_{\Delta'} / \left\langle \gamma_{i} \right\rangle} 
\left( \frac{1}{z - \gamma(\alpha_{i})} - \frac{1}{z - \gamma(\alpha_{-i})} \right) dz \ \ 
(1 \leq i \leq g-1), 
$$
and hence $\omega_{i}|_{y_{g} = 0}$ $(1 \leq i \leq g-1)$ are 
the universal abelian differentials of the first kind on ${\mathcal C}_{\Delta'}$ 
which give the basis of the space of relative abelian differentials of the first kind on ${\mathcal R}_{\Delta'}$ 
normalized for $\{ a_{i} \}_{1 \leq i \leq g-1}$. 
Furthermore, 
\begin{eqnarray*}
\omega_{g}|_{y_{g} = 0} 
& = & 
\sum_{\gamma \in \Gamma_{\Delta'}} 
\left( \frac{1}{z - \gamma(\alpha_{g})} - \frac{1}{z - \gamma(\alpha_{-g})} \right) dz 
\\ 
& = & 
\sum_{\gamma \in \Gamma_{\Delta'}} 
\left( \frac{d \gamma(z)}{\gamma(z) - \alpha_{g}} - \frac{d \gamma(z)}{\gamma(z) - \alpha_{-g}} \right), 
\end{eqnarray*}
and hence $\omega_{g}|_{y_{g} = 0}$ is the universal abelian differential of the third kind 
on ${\mathcal C}_{\Delta'}$ which gives the unique relative abelian differential of the third kind 
on ${\mathcal R}_{\Delta'}$ normalized for $\left( \{ a_{i} \}, \alpha_{g}, \alpha_{-g} \right)$. 
Therefore, the assertion (1) holds for $\left( {\mathcal R}_{\Delta}, \{ a_{i}, b_{i} \} \right)$. 

For an integer $k > 1$ and the unique tail $t$ of $\Delta$, 
by Theorem 3.4, the universal abelian differential of the second kind 
$$
\omega_{t, k} = 
\sum_{\gamma \in \Gamma_{\Delta}} \frac{d \gamma(z)}{(\gamma(z) - x_{t})^{k}} 
$$
on ${\mathcal C}_{\Delta}$ gives the unique relative abelian differential 
of the second kind with order $k$ on ${\mathcal R}_{\Delta}$ normalized for 
$\left( \{ a_{i} \}_{1 \leq i \leq g}, p_{t}, z - p_{t} \right)$. 
By Proposition 3.1 (3), 
$d \gamma(z)|_{y_{g}= 0} = 0$ if $\gamma \in \Gamma_{\Delta} - \Gamma_{\Delta'}$. 
Therefore, 
$$
\omega_{t, k}|_{y_{g} = 0} = 
\sum_{\gamma \in \Gamma_{\Delta'}} \frac{d \gamma(z)}{(\gamma(z) - x_{t})^{k}}, 
$$
and hence $\omega_{t, k}|_{y_{g} = 0}$ is the universal abelian differential of the second kind 
on ${\mathcal C}_{\Delta'}$ which gives the unique relative abelian differential of the second kind  with order $k$ on ${\mathcal R}_{\Delta'}$ 
normalized for $\left( \{ a_{i} \}_{1 \leq i \leq g-1}, p_{t}, z - x_{t} \right)$. 
Therefore, 
the assertion (2) holds for  $\left( {\mathcal R}_{\Delta}, \{ a_{i}, b_{i} \}, p_{t} \right)$, 
where $z - p_{t}$ can be replaced with general local coordinates at $p_{t}$. 

To prove the assertions generally, 
for integers $g, n \geq 0$ such that $2g - g + n > 0$, 
we consider the moduli space ${\mathcal H}_{g,n}$ 
of $(R, \{ a_{i}, b_{i} \}_{1 \leq i \leq g}, p_{1},..., p_{n})$, 
where $R$ are Riemann surfaces, $\{ a_{i}, b_{i} \}$ are symplectic bases of $R$ 
and $p_{1},..., p_{n}$ are points on $R$ which are different each other. 
Then ${\mathcal H}_{g,n}$ is a complex manifold as a covering space of 
the moduli space of Riemann surfaces of genus $g$ with $n$ marked points (cf. \cite{KnM, Kn}), 
and is seen to be connected by the Teichm\"{u}ller theory. 
Furthermore, 
the above $({\mathcal R}_{\Delta'}, \{ a_{i}, b_{i} \}_{1 \leq i \leq g-1}, p_{t_{1}}, p_{t_{2}}, z)$ 
make an open subset of ${\mathcal H}_{g-1,3}$, 
and $({\mathcal R}_{\Delta}, \{ a_{i}, b_{i} \}, z)$ give deformations of 
$\left( {\mathcal R}_{\Delta'}/(p_{t_{1}} = p_{t_{2}}), \{ a_{i}, b_{i} \}, z \right)$ 
for which the assertion (1) holds. 
Take a member $(R', \{ a'_{i}, b'_{i} \}, p_{1}, p_{2}, z)$ of ${\mathcal H}_{g-1,3}$ 
with deformation $({\mathcal R}, \{ a_{i}, b_{i} \}, z)$ of $\left( R'/(p_{1} = p_{2}), \{ a'_{i}, b'_{i} \}, z \right)$ 
by a complex parameter $y$ as above. 
Let $\{ \omega_{1},..., \omega_{g} \}$ be the basis of the space of relative abelian differentials of the first kind 
on ${\mathcal R}$ normalized for $\{ a_{i} \}$, 
$\{ \omega'_{1},..., \omega'_{g-1} \}$ be the basis of the space of abelian differentials of the first kind 
on $R'$ normalized for $\{ a'_{i} \}$, 
and $\omega'_{g}$ be the abelian differential of the third kind on $R'$ 
normalized for $(\{ a'_{i} \}, p_{1}, p_{2})$. 
Then for each $1 \leq i \leq g$, 
values at $z \in R'$ of 
$$
\lim_{y \rightarrow 0} \omega_{i} - \omega'_{i}
$$ 
are holomorphic on the moduli space ${\mathcal H}_{g-1,3}$ of $(R', \{ a'_{i}, b'_{i} \}, p_{1}, p_{2}, z)$ 
and are $0$ on its nonempty open subset 
$\left\{ ({\mathcal R}_{\Delta'}, \{ a_{i}, b_{i} \}, p_{t_{1}}, p_{t_{2}}, z) \right\}$. 
Therefore, by the identity theorem, 
these values are $0$ on the whole ${\mathcal H}_{g-1,3}$, 
and hence the assertion (1) holds generally. 

In order to show the assertion (2), for each positive integer $k$, 
we consider the moduli space ${\mathcal H}_{g-1, 2, k}$ of $(R', \{ a'_{i}, b'_{i} \}, p, u, z)$, 
where $(R', \{ a'_{i}, b'_{i} \}, p, z)$ are members of ${\mathcal H}_{g-1, 2}$ 
and $u$ are $k$-jets of local coordinates at $p$. 
Then ${\mathcal H}_{g-1, 2, k}$ is a connected complex manifold 
which contains $\left\{ ({\mathcal R}_{\Delta'}, \{ a_{i}, b_{i} \}, p_{t}, u, z) \right\}$ 
as its open subset, where $u$ are $k$-jets of local coordinates at $p_{t}$. 
Therefore, a similar argument as above shows the assertion (2). 
\ $\square$ 
\vspace{2ex}

\noindent
{\it 5.2. Reducible degeneration.} 
We consider abelian differentials on a general family of Riemann surfaces stably degenerating  to a reducible singular complex curve which is a union of two Riemann surfaces. 

Let ${\mathcal R}$ be a family of Riemann surfaces of genus $g > 0$ 
with symplectic bases $\{ a_{i}, b_{i} \}_{1 \leq i \leq g}$ and local coordinates $u$ at marked points $p$ 
which is defined over a small disk centered at $y = 0$ in ${\mathbb C}$. 
Take integers $g_{1}, g_{2} \geq 0$ satisfying $g_{1} + g_{2} = g$, 
and assume that under $y \rightarrow 0$, 
${\mathcal R}$ degenerates to a union of two Riemann surfaces $R_{1}, R_{2}$ 
of genus $g_{1}, g_{2}$ respectively by identifying $p_{1} \in R_{1}$ and $p_{2} \in R_{2}$, 
and $\{ a_{i}, b_{i} \}_{1 \leq i \leq g}$ becomes, as a homology basis, to the union of symplectic bases 
$$
\left\{ a_{1,i}, b_{1,i} \right\}_{1 \leq i \leq g_{1}}, \ 
\left\{ a_{2,i}, b_{2,i} \right\}_{g_{1} + 1 \leq i \leq g_{1} + g_{2}}
$$ 
of $R_{1}, R_{2}$ respectively. 
Let $\{ \omega_{1},..., \omega_{g} \}$ be the basis of the space of relative abelian differentials 
of the first kind on ${\mathcal R}$ normalized for $\{ a_{i} \}_{1 \leq i \leq g}$, 
and $\omega^{(n)}$ be the unique relative abelian differential of the second kind 
with order $n+1$ on ${\mathcal R}$ normalized for $(\{ a_{i} \}_{1 \leq i \leq g}, p, u)$. 
\vspace{2ex}

\noindent
{\bf Theorem 5.2.} 
\begin{it} 

{\rm (1)} 
Under $y \rightarrow 0$, 
$\{ \omega_{1},..., \omega_{g} \}$ becomes 
$$
\left\{ \omega_{1,1},..., \omega_{1,g_{1}}, 
\omega_{2,g_{1}+1},..., \omega_{2,g_{1}+g_{2}} \right\}, 
$$
where $\left\{ \omega_{1,1},..., \omega_{1,g_{1}} \right\}$ 
(resp. $\left\{ \omega_{2,g_{1}+1},..., \omega_{2,g_{1}+g_{2}} \right\}$) 
are the bases of the spaces of abelian differentials of the first kind on $R_{1}$ (resp. $R_{2}$) 
normalized for $\left\{ a_{1,i} \right\}_{1 \leq i \leq g_{1}}$ 
(resp. $\left\{ a_{2,i} \right\}_{g_{1}+1 \leq i \leq g_{1}+g_{2}}$). 

{\rm (2)} 
Assume that under $y \rightarrow 0$, 
$(p, u)$ becomes a point on $R_{1}$ with local coordinate which we denote by the same symbol. 
Then $\lim_{y \rightarrow 0} \omega^{(n)}$ is the unique abelian differential of the second kind 
with order $n+1$ on $R_{1}$ normalized for $\left( \left\{ a_{1,i} \right\}_{1 \leq i \leq g_{1}}, p, u \right)$. 
\end{it}
\vspace{2ex}

\noindent
{\it Proof.}
As is shown in the proof of Theorem 5.1, 
it is enough to prove the assertions (1) and (2) 
in the case when ${\mathcal R}$ is obtained from a generalized Tate curve as in 2.2 (P4). 
Let $\Delta_{i} = (V_{i}, E_{i}, T_{i})$ $(i = 1, 2)$ be two stable graphs 
such that ${\rm rank}_{\mathbb Z} H_{1}(\Delta_{i}, {\mathbb Z}) = g_{i}$, 
$V_{i} = \{ v_{i} \}$ and $T_{1} = \{ t, t_{1} \}$, $T_{2} = \{ t_{2} \}$, 
and $\Delta$ be a stable graph containing a nonloop oriented edge $e$ which is 
obtained from $\Delta_{1}, \Delta_{2}$ by connecting $t_{1}, t_{2}$ with $e$, 
where $v_{e} = v_{1}, v_{-e} = v_{2}$. 
Take generators $\gamma_{1,1},..., \gamma_{1,g_{1}}$ 
(resp. $\gamma_{2,g_{1}+1},..., \gamma_{2,g_{1}+g_{2}}$) of 
$\Gamma_{\Delta_{1}} \cong \pi_{1}(\Delta_{1}, v_{1})$ 
(resp. $\Gamma_{\Delta_{2}} \cong \pi_{1}(\Delta_{2}, v_{2})$), 
and symplectic bases $\{ a_{1,i}, b_{1,i} \}_{1 \leq i \leq g_{1}}$ 
(resp. $\{ a_{2,i}, b_{2,i} \}_{g_{1}+1 \leq i \leq g_{1}+g_{2}}$) 
of families ${\mathcal R}_{\Delta_{1}}$ (resp. ${\mathcal R}_{\Delta_{2}}$) 
of Riemann surfaces obtained from ${\mathcal C}_{\Delta_{1}}$ (resp. ${\mathcal C}_{\Delta_{2}}$) 
such that $b_{1,i}$ (resp. $b_{2,i}$) correspond to $\gamma_{1,i}$ (resp. $\gamma_{2,i}$). 
If $\phi_{e}$ is defined as in 2.2, then 
$$
\gamma_{i} = \left\{ \begin{array}{ll} 
\gamma_{1,i} & (1 \leq i \leq g_{1}), 
\\
\phi_{e} \cdot \gamma_{2,i} \cdot \phi_{e}^{-1} & (g_{1}+1 \leq i \leq g_{1}+g_{2}) 
\end{array} \right. 
$$
are generators of $\Gamma_{\Delta} \cong \pi_{1}(\Delta, v_{1})$. 
Therefore, by Theorem 3.4, the universal abelian differentials of the first kind 
$$
\omega_{i} = 
\sum_{\gamma \in \Gamma_{\Delta} / \left\langle \gamma_{i} \right\rangle} 
\left( \frac{1}{z - \gamma(\alpha_{i})} - \frac{1}{z - \gamma(\alpha_{-i})} \right) dz \ \ 
(1 \leq i \leq g) 
$$
on ${\mathcal C}_{\Delta}$ give the basis of the space of relative abelian differentials of the first kind 
on ${\mathcal R}_{\Delta}$ normalized for 
$$
\{ a_{1},..., a_{g} \} = \{ a_{1,1},..., a_{1,g_{1}}, a_{2,g_{1}+1},..., a_{2,g_{1}+g_{2}} \}. 
$$
By Proposition 3.1 (2), 
$$ 
\omega_{i}|_{y_{e} = 0} = \left\{ \begin{array}{l} 
\mbox{$\displaystyle 
\sum_{\gamma \in \Gamma_{\Delta_{1}} / \left\langle \gamma_{i} \right\rangle} 
\left( \frac{1}{z - \gamma(\alpha_{i})} - \frac{1}{z - \gamma(\alpha_{-i})} \right) dz$ 
on ${\mathcal C}_{\Delta_{1}}$ for $1 \leq i \leq g_{1}$,} 
\\ 
\mbox{$0$ on ${\mathcal C}_{\Delta_{2}}$ for $g_{1} + 1 \leq i \leq g_{1} + g_{2}$,} 
\end{array} \right. 
$$
and hence $\left\{ \omega_{i}|_{y_{e} = 0} \right\}_{1 \leq i \leq g_{1}}$ 
gives the basis of the space of relative abelian differentials of the first kind on ${\mathcal R}_{\Delta_{1}}$ 
normalized for $\{ a_{1,i} \}_{1 \leq i \leq g_{1}}$. 
Similarly, $\left\{ \omega_{i}|_{y_{e} = 0} \right\}_{g_{1}+1 \leq i \leq g_{1}+g_{2}}$ 
gives the basis of the space of relative abelian differentials of the first kind on ${\mathcal R}_{\Delta_{2}}$ 
normalized for $\{ a_{2,i} \}_{g_{1}+1 \leq i \leq g_{1}+g_{2}}$. 
For an integer $k > 1$ and the unique tail $t$ of $\Delta$, 
by Theorem 3.4, the universal abelian differential of the second kind 
$$
\omega_{t, k} = 
\sum_{\gamma \in \Gamma_{\Delta}} \frac{d \gamma(z)}{(\gamma(z) - x_{t})^{k}} 
$$
on ${\mathcal C}_{\Delta}$ gives the unique relative abelian differential of the second kind 
with order $k$ on ${\mathcal R}_{\Delta}$ normalized for 
$(\{ a_{i} \}_{1 \leq i \leq g}, p_{t}, z - p_{t})$. 
By Proposition 3.1 (3), 
$$ 
\omega_{t, k}|_{y_{e} = 0} = \left\{ \begin{array}{l} 
\mbox{$\displaystyle 
\sum_{\gamma \in \Gamma_{\Delta_{1}}} \frac{d \gamma(z)}{(\gamma(z) - x_{t})^{k}}$ 
on ${\mathcal C}_{\Delta_{1}}$,} 
\\ 
\mbox{$0$ on ${\mathcal C}_{\Delta_{2}}$,} 
\end{array} \right. 
$$
and hence $\omega_{t, k}|_{y_{e} = 0}$ gives the unique relative abelian differential 
of the second kind with order $k$ on ${\mathcal R}_{\Delta_{1}}$ normalized for 
$(\{ a_{i} \}_{1 \leq i \leq g-1}, p_{t}, z - x_{t})$. 
Therefore, the assertions (1) and (2) hold for $({\mathcal R}_{\Delta}, \{ a_{i}, b_{i} \}, p_{t})$, 
where $z - p_{t}$ is replaced with general local coordinates at $p_{t}$. 
This completes the proof. 
\ $\square$ 
\vspace{4ex}

\noindent
{\bf 6. Tau functions and quasi-periodic KP solutions}  
\vspace{2ex}

\noindent
{\it 6.1. Tau function and quasi-periodic KP solution.} 
We fix a Riemann surface $R$ of genus $g > 0$ 
with symplectic basis $\{ a_{i}, b_{i} \}_{1 \leq i \leq g}$, 
a point $p \in R$ and a local coordinate $u$ at $p$, 
and denote the whole data by 
$X = \left( R, \{ a_{i}, b_{i} \}_{1 \leq i \leq g}, p, u \right)$. 
Then there exists a unique basis $\{ \omega_{1},..., \omega_{g} \}$ of the space of abelian differentials 
of the first kind on $R$ normalized for $\{ a_{i} \}$, 
namely $\int_{a_{i}} \omega_{j} = 2 \pi \sqrt{-1} \delta_{ij}$ $(1 \leq i, j \leq g)$, 
and
$$
Z(X) = \left( Z_{i,j}(X) \right)_{1 \leq i, j \leq g} = 
\left( \frac{1}{2 \pi \sqrt{-1}} \int_{b_{i}} \omega_{j} \right)_{1 \leq i, j \leq g} 
$$ 
is the period matrix of $\left( R, \{ a_{i}, b_{i} \}_{1 \leq i \leq g} \right)$. 
Take $r_{j,m}(X) \in {\mathbb C}$ such that 
$$
\omega_{j} = \sum_{m=1}^{\infty} r_{j,m}(X) u^{m-1} du \ \ \mbox{at $p$}, 
$$
and put $r_{m}(X) = \left( r_{1,m}(X),..., r_{g,m}(X) \right)$. 
For each positive integer $n$, 
there exists a unique abelian differential $\omega^{(n)}$ of the second kind on $R$ satisfying 
\begin{itemize}

\item 
$\omega^{(n)}$ is holomorphic outside $p$, 

\item 
$\displaystyle \omega^{(n)} = 
\left( \frac{1}{u^{n+1}} + \sum_{m=1}^{\infty} \frac{q_{n,m}(X)}{n} u^{m-1} \right) du$ 
at $p$ for some $q_{n,m}(X) \in {\mathbb C}$, 

\item 
$\displaystyle \int_{a_{i}} \omega^{(n)} = 0$ for any $i = 1,..., g$. 

\end{itemize}
For $z = (z_{i})_{1 \leq i \leq g} \in {\mathbb C}^{g}$, 
we denote the Riemann theta function of $\left( R, \{ a_{i}, b_{i} \}_{1 \leq i \leq g} \right)$ by
$$
\Theta_{(R, \{ a_{i}, b_{i} \})} (z) = \sum_{v \in {\mathbb Z}^{g}} 
\left\{ \prod_{i,j = 1}^{g} \exp \left( \pi \sqrt{-1} Z_{i,j}(X) \right)^{v_{i} v_{j}} 
\prod_{i = 1}^{g} \exp(z_{i})^{v_{i}} \right\}, 
\eqno(6.1.1) 
$$
where $v = (v_{i})_{1 \leq i \leq g}$. 
Let $t_{m}$ $(m = 1, 2,...)$ be indeterminates, 
and put $t = (t_{1}, t_{2},...)$. 
Following \cite[Definition 5.5]{KNTY} (cf. \cite[(14)]{Al-GGR}, \cite{IMO}, \cite[(24)]{V}),  
for the addressed Riemann surface $X$ with line bundle of degree $0$ corresponding to 
$c = (c_{1},..., c_{g}) \in {\mathbb C}^{g}$, 
the associated {\it tau function} $\tau(t, X_{c})$ is defined as 
$$
\tau (t, X_{c}) = \exp \left( \frac{1}{2} \sum_{n, m = 1}^{\infty} q_{n,m}(X) t_{n} t_{m} \right) 
\cdot \Theta_{(R, \{ a_{i}, b_{i} \})} \left( c + \sum_{m = 1}^{\infty} r_{m}(X) t_{m} \right). 
\eqno(6.1.2) 
$$
Then $\tau(t, X_{c})$ is regarded as an element of 
${\mathbb C}[[t]] = {\mathbb C} [[t_{1}, t_{2},...]]$ by using the expression 
$$
\prod_{i = 1}^{g} \exp \left( \sum_{m = 1}^{\infty} r_{i,m}(X) t_{m} \right)^{v_{i}} = 
\sum_{n = 0}^{\infty} \frac{1}{n!} 
\left( \sum_{m = 1}^{\infty} \left( \sum_{i = 1}^{g} v_{i} \cdot r_{i,m}(X) \right) t_{m} \right)^{n}. 
$$
If $\Theta_{(R, \{ a_{i}, b_{i} \})}(c) \neq 0$, 
then we put
$$
\frac{\tau(t - [\alpha], X_{c})}{\tau(t, X_{c})} = 
1 + \sum_{k = 1}^{\infty} w_{k} (t, X_{c}) \alpha^{k}, 
\eqno(6.1.3) 
$$
where $[\alpha] = \left( \alpha, \alpha^{2}/2, \alpha^{3}/3,... \right)$, 
and define two micro-differential operators
$$
W (t, X_{c}) = 1 + \sum_{k = 1}^{\infty} w_{k} (t, X_{c}) \partial_{x}^{-k}
\eqno(6.1.4) 
$$
and
$$
L (t, X_{c}) = W (t + x, X_{c}) \cdot \partial_{x} \cdot W (t + x, X_{c})^{-1}
\eqno(6.1.5) 
$$ 
with coefficients in ${\mathbb C}[[x, t]]$, 
where $t + x = (t_{1} + x, t_{2}, t_{3},... )$. 
Then it is known (cf. \cite{Kr, SeW}) that $L(t, X_{c})$ satisfies the KP hierarchy
$$
\frac{\partial L}{\partial t_{n}} = \left[ (L^{n})_{+}, L \right] \ \ (n = 1, 2,...). 
$$
In particular,
$$
u_{1} (x, t_{2}, t_{3}) = \frac{\partial^{2}}{\partial x^{2}} 
\log \Theta_{(R, \{ a_{i}, b_{i} \})} \left( c + x r_{1}(X) + t_{2} r_{2}(X) + t_{3} r_{3}(X) \right) +  q_{1,1}(X) 
$$
satisfies the KP equation
$$
\frac{3}{4} \frac{\partial^{2} u_{1}}{\partial t_{2}^{2}} - 
\frac{\partial}{\partial x} 
\left( \frac{\partial u_{1}}{\partial t_{3}} - \frac{1}{4} \frac{\partial^{3} u_{1}}{\partial x^{3}} - 
3 u_{1} \frac{\partial u_{1}}{\partial x} \right) = 0. 
$$ 
Since $L(t, X_{c})$ and $u_{1}(x, t_{2}, t_{3})$ are expressed by theta functions 
with quasi-periodicity, 
they are called {\it quasi-periodic KP solutions}. 
\vspace{2ex}

\noindent
{\it 6.2. Universal tau function.} 
Let $\Delta = (V, E, T)$ be a stable graph of $(g, 1)$-type with $T = \{ t_{0} \}$, 
and fix generators $\rho_{1},..., \rho_{g}$ of $\pi_{1}(\Delta)$. 
Then we denote by ${\mathcal X}_{\Delta}$ the family ${\mathcal R}_{\Delta}$ 
of Riemann surfaces associated with ${\mathcal C}_{\Delta}$ 
with symplectic basis $\{ a_{i}, b_{i} \}_{1 \leq i \leq g}$, marked point $x_{t_{0}} \in P_{v_{t_{0}}}$ 
and local coordinate $z - x_{t_{0}}$, 
where $b_{i}$ correspond to $[\rho_{i}] \in H_{1}(\Delta, {\mathbb Z})$ $(i = 1,..., g)$.  
\vspace{2ex}

\noindent
{\bf Theorem 6.1.} 
\begin{it} 

{\rm (1)} 
Denote by ${\mathbb Q} \left[ \gamma_{1}^{\pm 1},..., \gamma_{g}^{\pm 1} \right]$ 
the ${\mathbb Q}$-algebra of Laurent polynomials in variables $\gamma_{1},..., \gamma_{g}$. 
Then there exists the universal tau function 
$\tau \left( t, ({\mathcal C}_{\Delta})_{\gamma} \right)$ as an element of 
$$
{\mathcal A}_{\Delta} = 
A_{\Delta} \left[ \sqrt{P_{ii}} \ (i = 1,..., g) \right] \widehat{\otimes}_{\mathbb Z} 
{\mathbb Q} \left[ \gamma_{1}^{\pm 1},..., \gamma_{g}^{\pm 1} \right] [[t]] 
$$
which becomes the tau functions of members of $({\mathcal X}_{\Delta})_{c}$, 
where $\exp(c_{i}) = \gamma_{i}$. 

{\rm (2)} 
The universal tau function $\tau \left( t, ({\mathcal C}_{\Delta})_{\gamma} \right)$ 
gives rise to a KP solution as in 6.1. 

{\rm (3)} 
By 2.2 (P2),  $\tau \left( t, ({\mathcal C}_{\Delta})_{\gamma} \right)$ gives tau functions 
of all marked Mumford curves with dual graph $\Delta$ 
over nonarchimedean complete valuation fields of characteristic $0$, 
and these tau functions give rise to KP solutions as in 6.1. 

\end{it} 
\vspace{2ex}

\noindent
{\it Proof.} 
The assertion (1) follows from Theorems 3.4 and 3.6. 
By \cite[Proposition 3.11]{I2}, 
$P_{ii}$ $(i = 1,..., g)$ belong to the ideal of $A_{\Delta}$ generated by $y_{e}$ $(e \in E)$, 
and $\tau \left( t, ({\mathcal C}_{\Delta})_{\gamma} \right) - 1$ belongs to the ideal of 
${\mathcal A}_{\Delta}$ generated by $\sqrt{P_{ii}}$ $(i = 1,..., g)$, 
and hence $\tau \left( t, ({\mathcal C}_{\Delta})_{\gamma} \right)$ 
defines a micro-differential operator $L \left( t, ({\mathcal C}_{\Delta})_{\gamma} \right)$ 
with coefficients in ${\mathcal A}_{\Delta}[[x]]$. 
Then the assertion (2) follows from (1), 
and (3) follows from (2). 
\ $\square$ 
\vspace{2ex} 

\noindent
{\it Remark 6.2.} 
By the results in Section 4, 
one can see that 
$\tau \left( t, ({\mathcal C}_{\Delta})_{\gamma} \right)$ 
($\Delta$: stable graphs of $(g, 1)$-type) gives rise to tau functions 
over the family of all degenerating Riemann surfaces, 
and that of all Mumford curves. 
\vspace{4ex}

\noindent
{\bf 7. Behavior of tau functions}  
\vspace{2ex}

\noindent
{\it 7.1. Variation of tau functions.} 
Let the notation be as in Section 6, 
and put 
$$
c = 2 \pi \sqrt{-1} (\alpha + Z(X) \beta), 
$$
where 
$\alpha = (\alpha_{i})_{1 \leq i \leq g} \in {\mathbb C}^{g}, 
\beta = (\beta_{i})_{1 \leq i \leq g} \in {\mathbb R}^{g}$. 
Then we study the asymptotic behavior of the tau function (6.1.2) 
under degenerations of $X$ as in Section 6, 
and define the associated modified tau function as its regularized limit. 
\vspace{2ex}

\noindent
{\it 7.2. Irreducible degeneration.} 
Let the notation be as in 5.1, especially 
$$
{\mathcal X} = \left( {\mathcal R}, \{ a_{i}, b_{i} \}_{1 \leq i \leq g}, p, u \right), \ \ 
X' = \left( R', \{ a'_{i}, b'_{i} \}_{1 \leq i \leq g-1}, p, u \right), 
$$
and ${\mathcal R}$ degenerates to $R'$ by $y \rightarrow 0$. 
Then we study the relationship between the the tau functions of ${\mathcal X}$ and $X'$. 
\vspace{2ex}

\noindent
{\bf Theorem 7.1.} 
\begin{it} 
Assume that $\beta_{g} \not\in {\mathbb Z} + 1/2$, 
and take a unique integer $\overline{\beta}$ such that 
$\left| \overline{\beta} + \beta_{g} \right| < 1/2$. 

{\rm (1)} 
Under $y \rightarrow 0$, 
$$
\exp \left( - \pi \sqrt{-1} \left( \overline{\beta} + 2 \beta_{g} \right) \overline{\beta} 
Z_{g,g}({\mathcal X}) \right) \cdot 
\Theta_{({\mathcal R}, \{ a_{i}, b_{i} \})} 
\left( c + \sum_{m=1}^{\infty} r_{m}({\mathcal X}) t_{m}\right)
$$
tends to 
$$
\exp \left( 2 \pi \sqrt{-1} \left( \alpha_{g} + 
\sum_{i=1}^{g-1} \beta_{i} \overline{Z}_{i,g} \right) + 
\sum_{m=1}^{\infty} \overline{r}_{g,m} t_{m} \right)^{\overline{\beta}} \cdot 
\Theta_{(R', \{ a'_{i}, b'_{i} \})} 
\left( \overline{c} + \sum_{m=1}^{\infty} \overline{r}_{m} t_{m} \right), 
$$
where 
\begin{eqnarray*}
\overline{Z}_{i,g} 
& = & 
\lim_{y \rightarrow 0} Z_{i,g}({\mathcal X}) 
\ \ (1 \leq i \leq g-1), 
\\ 
\overline{c} 
& = & 
\left( \alpha_{i} + \overline{Z}_{i,g} (\overline{\beta} + \beta_{g}) \right)_{1 \leq i \leq g-1}, 
\\ 
\overline{r}_{m} 
& = & 
\left( r_{i,m} \left( X' \right) \right)_{1 \leq i \leq g-1}, 
\end{eqnarray*}
and $\sum_{m=1}^{\infty} \overline{r}_{g,m} u^{m-1} du$ is the expansion of 
$\omega'$ given in Theorem 5.1 (1). 

{\rm (2)} 
Under $y \rightarrow 0$, 
$$
\exp \left( - \pi \sqrt{-1} \left( \overline{\beta} + 2 \beta_{g} \right) \overline{\beta} 
Z_{g,g}({\mathcal X}) \right) \cdot 
\tau \left( t, {\mathcal X}_{c} \right) 
$$
tends to 
$$
\exp \left( 2 \pi \sqrt{-1} \left( \alpha_{g} + 
\sum_{i=1}^{g-1} \beta_{i} \overline{Z}_{i,g} \right) + 
\sum_{m=1}^{\infty} \overline{r}_{g,m} t_{m} \right)^{\overline{\beta}} \cdot 
\tau \left( t, X'_{\overline{c}} \right) 
$$
which we call the associated modified tau function. 
\end{it}
\vspace{2ex}

\noindent
{\it Proof.} 
By (6.1.1), the assertion (1) follows from the formula 
\begin{eqnarray*} 
& & 
\Theta \left( c + \sum_{m=1}^{\infty} r_{m} t_{m} \right)
\\ 
& = & 
\sum_{v \in {\mathbb Z}^{g}} 
\left\{ \prod_{i, j=1}^{g} \exp \left( \pi \sqrt{-1} v_{i} Z_{i,j} v_{j} \right) 
\prod_{i=1}^{g} \exp \left( 2 \pi \sqrt{-1} \alpha_{i} \right)^{v_{i}} 
\prod_{i, j=1}^{g} \exp \left( 2 \pi \sqrt{-1} v_{i} Z_{i,j} \beta_{j} \right) 
\right. 
\\ 
& & 
\left. 
\times \prod_{i=1}^{g} \exp \left( \sum_{m=1}^{\infty} r_{i,m} t_{m} \right)^{v_{i}} 
\right\} 
\\ 
& = & 
\prod_{i, j=1}^{g} \exp \left( - \pi \sqrt{-1} \beta_{i} Z_{i,j} \beta_{j} \right) 
\sum_{v \in {\mathbb Z}^{g}} 
\left\{ \prod_{i, j=1}^{g} \exp \left( \pi \sqrt{-1} (v_{i} + \beta_{i}) Z_{i,j} (v_{j} + \beta_{j}) \right) 
\right. 
\\ 
& & 
\left. 
\times 
\prod_{i=1}^{g} \exp \left( 2 \pi \sqrt{-1} \alpha_{i} \right)^{v_{i}} 
\sum_{n=0}^{\infty} \frac{1}{n!} 
\left( \sum_{m=1}^{\infty} \left( \sum_{i=1}^{g} v_{i} r_{i,m} \right) t_{m} \right)^{n} 
\right\}, 
\end{eqnarray*}
and the assertion (2) follows from (1). 
\ $\square$ 
\vspace{2ex}

\noindent
{\bf Theorem 7.2.} 
\begin{it} 
Assume that $\beta_{g} \in {\mathbb Z} + 1/2$, 
and put $\overline{\beta}_{1} = - \beta_{g} - 1/2, \overline{\beta}_{2} = - \beta_{g} + 1/2$ 
which are the nearest integers to $\beta_{g}$. 

{\rm (1)} 
Under $y \rightarrow 0$, 
$$
\exp \left( \pi \sqrt{-1} \left( {\beta}_{g}^{2} - 1/4 \right) Z_{g,g}({\mathcal X}) \right) \cdot 
\Theta_{({\mathcal R}, \{ a_{i}, b_{i} \})} 
\left( c + \sum_{m=1}^{\infty} r_{m}(X) t_{m} \right)
$$
tends to 
\begin{eqnarray*}
& & 
\exp \left( 2 \pi \sqrt{-1} \left( \alpha_{g} + 
\sum_{i=1}^{g-1} \beta_{i} \overline{Z}_{i,g} \right) + 
\sum_{m=1}^{\infty} \overline{r}_{g,m} t_{m} \right)^{\overline{\beta}_{1}} 
\\ 
& & 
\times \ 
\Theta_{(R', \{ a'_{i}, b'_{i} \})} 
\left( \overline{c}_{1} + \sum_{m=1}^{\infty} \overline{r}_{m} t_{m} \right) 
\\ 
& + & 
\exp \left( 2 \pi \sqrt{-1} \left( \alpha_{g} + 
\sum_{i=1}^{g-1} \beta_{i} \overline{Z}_{i,g} \right) + 
\sum_{m=1}^{\infty} \overline{r}_{g,m} t_{m} \right)^{\overline{\beta}_{2}} 
\\ 
& & 
\times \ 
\Theta_{(R', \{ a'_{i}, b'_{i} \})} 
\left( \overline{c}_{2} + \sum_{m=1}^{\infty} \overline{r}_{m} t_{m} \right), 
\end{eqnarray*}
where 
\begin{eqnarray*}
\overline{Z}_{i,g} 
& = & 
\lim_{y \rightarrow 0} Z_{i,g}({\mathcal X}) 
\ \ (1 \leq i \leq g-1), 
\\ 
\overline{c}_{k} 
& = & 
\left( \alpha_{i} + \overline{Z}_{i,g} (\overline{\beta}_{k} + \beta_{g}) \right)_{1 \leq i \leq g-1} 
\ \ (k = 1, 2), 
\\ 
\overline{r}_{m} 
& = & 
\left( r_{i,m} \left( X' \right) \right)_{1 \leq i \leq g-1}, 
\end{eqnarray*}
and $\sum_{m=1}^{\infty} \overline{r}_{g,m} u^{m-1} du$ is the expansion of 
$\omega'$ given in Theorem 5.1 (1). 

{\rm (2)} 
Under $y \rightarrow 0$, 
$$
\exp \left( \pi \sqrt{-1} \left( {\beta}_{g}^{2} - 1/4 \right) Z_{g,g}({\mathcal X}) \right) 
\cdot \tau \left( t, {\mathcal X}_{c} \right) 
$$
tends to 
\begin{eqnarray*}
& & 
\exp \left( 2 \pi \sqrt{-1} \left( \alpha_{g} + 
\sum_{i=1}^{g-1} \beta_{i} \overline{Z}_{i,g} \right) + 
\sum_{m=1}^{\infty} \overline{r}_{g,m} t_{m} \right)^{\overline{\beta}_{1}} 
\cdot 
\tau \left( t, X'_{\overline{c}_{1}} \right) 
\\ 
& + & 
\exp \left( 2 \pi \sqrt{-1} \left( \alpha_{g} + 
\sum_{i=1}^{g-1} \beta_{i} \overline{Z}_{i,g} \right) + 
\sum_{m=1}^{\infty} \overline{r}_{g,m} t_{m} \right)^{\overline{\beta}_{2}} 
\cdot 
\tau \left( t, X'_{\overline{c}_{2}} \right) 
\end{eqnarray*}
which we call the associated modified tau function. 
\end{it} 
\vspace{2ex}

\noindent
{\it Proof.} 
One can prove the assertion as in Theorem 7.1. 
\ $\square$ 
\vspace{2ex}

\noindent
{\it 7.3. Reducible degeneration.} 
Let the notation and assumption be as in 5.2, 
especially  
$$ 
{\mathcal X} = \left( {\mathcal R}, \{ a_{i}, b_{i} \}_{1 \leq i \leq g}, p, u \right), \ \ 
X_{1} = \left( R_{1}, \left\{ a_{1,i}, b_{1,i} \right\}_{1 \leq i \leq g_{1}}, p, u \right), 
$$
and ${\mathcal R}$ degenerates to the union of $R_{1}, R_{2}$ by $y \rightarrow 0$. 
\vspace{2ex}

\noindent
{\bf Theorem 7.3.} 
\begin{it} 
Let $c = (c_{1},..., c_{g})$ be an element of ${\mathbb C}^{g}$. 

{\rm (1)} 
Under $y \rightarrow 0$, 
$\Theta_{({\mathcal R}, \{ a_{i}, b_{i} \})} 
\left( c + \sum_{m = 1}^{\infty} r_{m}({\mathcal X}) t_{m} \right)$ tends to 
$$
\Theta_{(R_{1}, \{ a_{1,i}, b_{1,i} \})} 
\left( c' + \sum_{m = 1}^{\infty} r_{m}(X_{1}) t_{m} \right) 
\cdot \Theta_{(R_{2}, \{ a_{2,i}, b_{2,i} \})} \left( c'' \right), 
$$ 
where $c' = (c_{1},..., c_{g_{1}})$ and $c'' = (c_{g_{1}+1},..., c_{g_{1}+g_{2}})$. 

{\rm (2)} 
Under $y \rightarrow 0$, 
$\tau \left( t, {\mathcal X}_{c} \right)$ tends to 
$$
\tau \left( t, (X_{1})_{c'} \right) \cdot 
\Theta_{(R_{2}, \{ a_{2,i}, b_{2,i} \})} \left( c'' \right) 
$$ 
which we call the associated modified tau function. 
\end{it}
\vspace{2ex}

\noindent
{\it Proof.}
By Theorem 5.2,  
$$
\lim_{y \rightarrow 0} Z_{i,j} = \left\{ \begin{array}{ll} 
0 & (1 \leq i \leq g_{1}, g_{1} + 1 \leq j \leq g_{1} + g_{2}), 
\\ 
0 & (g_{1} + 1 \leq i \leq g_{1} + g_{2}, 1 \leq j \leq g_{1}), 
\\ 
\mbox{$(i, j)$ period of $\left( R_{1}, \{ a_{1,i}, b_{1,i} \} \right)$} 
& (1 \leq i, j \leq g_{1}), 
\\
\mbox{$(i, j)$ period of $\left( R_{2}, \{ a_{2,i}, b_{2,i} \} \right)$} 
& (g_{1} + 1 \leq i, j \leq g_{1} + g_{2}). 
\end{array} \right. 
$$
Therefore, the assertions (1) and (2) follow from Theorem 5.2. 
\ $\square$  
\vspace{2ex}

\noindent
{\it 7.4. KP solutions and generalized soliton.} 
By using Theorems 7.1, 7.2 and 7.3 repeatedly, 
for degenerations to stable complex curves with suitably chosen symplectic bases, 
we can express the associated modified tau functions as mixtures of quasi-periodic solutions and solitons. 
Then we have the following assertion. 
\vspace{2ex}

\noindent
{\bf Theorem 7.4.} 
\begin{it} 
The modified tau functions give rise to KP solutions as in 6.1. 
\end{it}
\vspace{2ex}

\noindent 
{\it Proof.} 
We only prove the assertion in the case of Theorem 7.1 since other cases can be shown similarly. 
By (6.1.3), 
$$
\tau \left( t, {\mathcal X}_{c} \right), \ \ 
\exp \left( - \pi \sqrt{-1} \left( \overline{\beta} + 2 \beta_{g} \right) \overline{\beta} 
Z_{g,g}({\mathcal X}) \right) \cdot 
\tau \left( t, {\mathcal X}_{c} \right) 
$$
gives the same micro-differential operators in (6.1.4) and (6.1.5), 
and hence their limits under $y \rightarrow 0$ is associated with the modified tau function, 
and gives rise to a KP solution. 
\ $\square$ 
\vspace{2ex}

In particular,
we have solutions expressed by rational functions containing soliton solutions as follows. 
Let $\Delta_{0} = (V_{0}, E_{0}, T_{0})$ be the stable graph of $(g, 1)$-type 
consisting of one vertex $v_{0}$, 
and identify $E_{0}$, $T_{0}$ with $\{ 1,..., g \}$, $\{ t_{0} \}$ respectively. 
Then under $y_{i} \rightarrow 0$ $(i = 1,..., g)$, 
${\mathcal X}_{\Delta_{0}}$ degenerates to to the irreducible complex curve 
${\mathbb P}^{1}/(x_{i} = x_{-i})$ which is obtained from ${\mathbb P}^{1} = P_{v_{0}}$ 
identifying $2g$ points $x_{\pm i}$ $(i = 1,..., g)$ 
in pairs $x_{i}, x_{-i}$ with local coordinate $u = z - x_{t_{0}}$ at $x_{t_{0}}$. 
\vspace{2ex}

\noindent
{\bf Theorem 7.5.} 
\begin{it} 
Let $\beta = \left( n_{i} - 1/2 \right)_{1 \leq i \leq g}$ be an element of 
$\left( {\mathbb Z} + 1/2 \right)^{g}$, 
and put 
$$
c = 2 \pi \sqrt{-1} \left( \alpha + Z \left( {\mathcal X}_{\Delta_{0}} \right) \beta \right)
$$ 
as in 7.1, 
Then there exists $\alpha' = (\alpha'_{i})_{1 \leq i \leq g} \in {\mathbb C}^{g}$ such that 
the associated modified tau function 
$$
\lim_{y_{1},..., y_{g} \rightarrow 0} 
\left( \prod_{i=1}^{g} \exp 
\left( \pi \sqrt{-1} (n_{i}^{2} - n_{i}) Z_{ii} \left( {\mathcal X}_{\Delta_{0}} \right) \right) 
\cdot \tau \left( t, \left( {\mathcal X}_{\Delta_{0}} \right)_{c} \right) \right) 
$$
is expressed as  
\begin{eqnarray*} 
& & 
\sum_{u \in \{ 0, 1 \}^{g}} \left\{ 
\prod_{1 \leq i < j \leq g} 
\left( \frac{(x_{i} - x_{j})(x_{-i} - x_{-j})}{(x_{-i} - x_{j})(x_{i} - x_{-j})} 
\right)^{(u_{i} - n_{i})(u_{j} - n_{j})} 
\right. 
\\ 
& &
\left. \times 
\prod_{i=1}^{g} 
\exp \left( \alpha'_{i} - \sum_{m=1}^{\infty} 
\left( \frac{1}{(x_{i} - x_{t_{0}})^{m}} - \frac{1}{(x_{-i} - x_{t_{0}})^{m}} \right)t_{m} \right)^{u_{i} - n_{i}} 
\right\} 
\end{eqnarray*}
which gives rise to KP solutions as in 6.1. 
When $n_{i} = 0$ $(i = 1,..., g)$, 
these solutions are the soliton solutions. 
\end{it}
\vspace{2ex}

\noindent 
{\it Proof.} 
By using Theorem 7.2 repeatedly, 
the modified tau function becomes the sum over $u \in \{ 0, 1 \}^{g}$ of 
$$
\prod_{1 \leq i < j \leq g} 
\exp \left( 2 \pi \sqrt{-1} (u_{i} - n_{i}) \overline{Z}_{i,j} (u_{j} - n_{j}) \right) 
\cdot 
\prod_{i=1}^{g} \exp \left( (u_{i} - n_{i}) 
\left( \alpha'_{i} + \sum_{m=1}^{\infty} \overline{r}_{i,m} t_{m} \right) \right). 
$$
Denote by $P_{ij}$ $(1 \leq i, j \leq g)$ the multiplicative periods 
of ${\mathcal C}_{\Delta_{0}}$ as in 3.3. 
By Theorem 3.6 (1), if $i \neq j$, then 
$$
P_{ij} = \frac{(x_{i} - x_{j})(x_{-i} - x_{-j})}{(x_{-i} - x_{j})(x_{i} - x_{-j})} \cdot 
\prod_{\gamma} \left( 1 + \frac{(x_{i} - x_{-i})(\gamma(x_{j}) - \gamma(x_{-j}))}
{(x_{-i} - \gamma(x_{j}))(x_{i} - \gamma(x_{-j}))} \right), 
$$
where $\gamma$ runs through 
$\langle \gamma_{i} \rangle \backslash 
\left( \Gamma_{\Delta_{0}} - \langle \gamma_{i} \rangle \langle \gamma_{i} \rangle \right) 
/ \langle \gamma_{j} \rangle$, 
and hence by Proposition 3.1 (2), 
$$
\exp \left( 2 \pi \sqrt{-1} \overline{Z}_{i,j} \right) 
= \lim_{y_{1},..., y_{g} \rightarrow 0} P_{ij} 
= \frac{(x_{i} - x_{j})(x_{-i} - x_{-j})}{(x_{-i} - x_{j})(x_{i} - x_{-j})}. 
$$
Furthermore, by Proposition 3.1 (2), 
the abelian differential $\omega_{i}$ $(i = 1,..., g)$ on ${\mathcal C}_{\Delta_{0}}$ 
given in 3.1 satisfy 
$$
\lim_{y_{1},..., y_{g} \rightarrow 0} \omega_{i} 
= - \left( \frac{1}{(x_{i} - x_{t_{0}}) - u} - \frac{1}{(x_{-i} - x_{t_{0}}) - u} \right) du, 
$$
and hence 
$$
\overline{r}_{i,m} = 
- \left( \frac{1}{(x_{i} - x_{t_{0}})^{m}} - \frac{1}{(x_{-i} - x_{t_{0}})^{m}} \right). 
$$
Then we have the above expression of the modified tau function. 
The remaining assertion follows from Theorem 7.4 and results in \cite{Go, Mu3}. 
\ $\square$ 
\vspace{4ex}

\noindent
{\bf 8. Family of {\texttt M}-curves and real tau functions}  
\vspace{2ex}

\noindent
{\it 8.1. Family of {\texttt M}-curves.} 
Let $\Delta_{0} = (V_{0}, E_{0}, T_{0})$ be the stable graph of $(g, n)$-type consisting of 
one vertex $v_{0}$ and $g$ oriented loops $\rho_{1},..., \rho_{g}$, 
and identify $E_{0}$ with $\{ 1,..., g \}$. 
Take $x_{\pm i}$ $(i = 1,..., g)$, $x_{t}$ $(t \in T_{0})$ 
as real numbers distinct to each other satisfying 
$$
x_{-1} < x_{1} < x_{-2} < x_{2} < \cdots < x_{-g} < x_{g}
$$
which are regarded as cyclically ordered real points on $P_{v_{0}} = {\mathbb P}^{1}$, 
$y_{i}$ $(i = 1,..., g)$ as sufficiently small positive real numbers. 
Then by the relation (2.2.1), 
the associated Schottky uniformized Riemann surfaces are {\texttt M}-curves, 
namely real algebraic curves with maximal number of components. 
Furthermore, 
under comparing ${\mathcal C}_{\Delta}$ with ${\mathcal C}_{\Delta'}$ as in Theorem 4.1, 
if the parameters $x_{h}$ $(h \in \pm E \cup T)$, $y_{e}$ $(e \in E)$ 
for ${\mathcal C}_{\Delta}$ are given as real numbers associated with {\texttt M}-curves, 
then we can construct {\texttt M}-curves from ${\mathcal C}_{\Delta'}$ 
by taking the parameters $s_{e}$ $(e \in E')$ as real numbers such that
the ratios in Theorem 4.1 are positive. 
Therefore, by applying the construction of the universal Mumford curve, 
we have a family ${\texttt M}_{g,n}$ of $n$-marked {\texttt M}-curves of genus $g$. 
\vspace{2ex}

\noindent
{\it 8.2. Real tau functions.} 
\vspace{2ex}

\noindent
{\bf Theorem 8.1.} 
\begin{it}
Let $\Delta_{0} = (V_{0}, E_{0}, T_{0})$ be the stable graph of $(g, 1)$-type. 
Then for each $\gamma \in ({\mathbb R}^{\times})^{g}$, 
the universal tau function $\tau \left( t, ({\mathcal C}_{\Delta_{0}})_{\gamma} \right)$ 
given in Theorem 6.1 (1) gives rise to tau functions defined for 
all members in ${\texttt M}_{g,1}$ which are real, 
namely real-valued for real variables $t_{m}$ $(m = 1, 2,...)$. 
\end{it}
\vspace{2ex}

\noindent
{\it Proof.} 
For the above {\texttt M}-curves defined by $\Delta_{0}$, 
Theorem 3.6 implies that 
$$
\exp \left( 2 \pi \sqrt{-1} Z_{i,j} \right) \in {\mathbb R}^{\times}, \ \ 
0 < \exp \left( \pi \sqrt{-1} Z_{i,i} \right) = 
\exp \left( \frac{1}{2} \int_{[\rho_{i}]} \omega_{i} \right) < 1 
$$
since $\int_{[\rho_{i}]}$ can be taken as the integral along a real line in 
${\mathcal X}_{\Delta_{0}}$. 
Therefore, the associated tau functions are real, 
and hence they are extended to real functions defined for all members in ${\texttt M}_{g, 1}$. 
\ $\square$ 
\vspace{2ex} 
 
\noindent
{\it Remark 8.2.} 
Using results in Section 7, especially Theorem 7.5, 
one can similarly construct real quasi-periodic KP solutions tending to soliton solutions 
studied in \cite{AG} which are also obtained from theta functions of tropical curves 
as their tropical limits in \cite{AFMS}. 
\vspace{4ex}

\noindent 
{\bf Acknowledgments} 
\vspace{2ex}

This work is partially supported by the JSPS Grant-in-Aid for 
Scientific Research No. 20K03516.

\renewcommand{\refname}{\normalsize{\bf References}}

\end{document}